\documentclass[11pt,a4paper,geompsfi]{article}
\usepackage{amsfonts}

\usepackage[all]{xy}
\usepackage[latin1]{inputenc}        
\usepackage[dvips]{graphics}
\usepackage[dvips]{graphicx}
\usepackage{amssymb}
\usepackage{amsmath}
\usepackage{amstext}
\usepackage{amsbsy}
\usepackage{amsopn}
\usepackage{amsthm}
\usepackage{epsfig,rotating}
\usepackage{amscd}
\usepackage{amsxtra}
\usepackage{mathrsfs}
\usepackage{upref}

\let\small=\footnotesize
\def\longpage#1{\newdimen\addstr  \addstr=#1\baselineskip
  \advance\topmargin-2\baselineskip \advance\textheight4\baselineskip
  \advance\topmargin-0.5\addstr     \advance\textheight\addstr
  \advance\oddsidemargin-0.5\addstr \advance\textwidth\addstr}
\longpage{0.53}

       \def\R{\mathbb R}        
              \def\F{\Phi}      
                  \def\f{\varphi}       
          \def\:{\colon\,}        
\def\D{\partial}               
\def\s{\sigma}                 \def\g{\gamma} 
\def\e{\varepsilon}

\def\P{{\cal P}}               
               
\def\comp{{\small{\circ}}}     \def\a{\alpha}
                 
\def\fin{\hfill{$\square$}}

               \def\A{{\mathcal A}}  
          
              \def\setminus{\smallsetminus}
             \def\ra{\rightarrow}
      \def\tilde{\widetilde}
          \def\bar{\overline}
\def\sF{{\scriptscriptstyle F}}  \def\sP{{\scriptscriptstyle P}}
\def\sD{{\scriptscriptstyle D}} \def\sT{{\scriptscriptstyle T}}
\def\sC{{\scriptscriptstyle C}}
\def\RP{\mathbb R\mathrm P}    
\def\:{{\rm :}}

\author{by Ricardo \sc Uribe-Vargas\footnote{Partially supported by EU
Centre of Excellence, IMPAN-Banach Centre, ICA1-CT-2000-70024.} 
\\ {\small Coll\`ege de France, 3 rue d'Ulm, 
75005 Paris.}\\ 
{\small uribe@math.jussieu.fr \ \ \ 
www.math.jussieu.fr/{\small $\sim$}uribe/ }}

\date\empty                     
\title{A New Projective Invariant for Swallowtails and Godrons 
(Cusps of Gauss), and Global Theorems on the 
Flecnodal Curve}
      
\begin{document}

\numberwithin{equation}{section}                
\theoremstyle{plain}
\newtheorem{theorem}{\bf Theorem}
\newtheorem*{Theorem}{\bf Theorem}
\newtheorem*{lrtheorem}{\bf Left-Right-Theorem}
\newtheorem*{Euler}{\bf Euler's Criterion}
\newtheorem{lemma}{\bf Lemma}
\newtheorem*{slemma}{\bf Separating $2$-jet Lemma}
\newtheorem*{qlemma}{\bf $q$-Contour Lemma}
\newtheorem*{lemma*}{\bf Lemma}
\newtheorem*{lemma0}{\bf Lemma 0}
\newtheorem*{lemmaa}{\bf Lemma A}
\newtheorem*{lemmab}{\bf Lemma B}
\newtheorem*{lemmac}{\bf Lemma C}
\newtheorem*{lemmad}{\bf Lemma D}
\newtheorem{claim}{\bf Claim}
\newtheorem*{fitheorem}{\bf Inflection-Fold Theorem}
\newtheorem{proposition}{\bf Proposition}
\newtheorem*{Proposition}{\bf Proposition}
\newtheorem{corollary}{\bf Corollary}
\newtheorem*{corollary3}{\bf Corollary {\rm (of Theorem \ref{General})}}

\theoremstyle{definition}
\newtheorem{definition}{\bf Definition}
\newtheorem*{definition*}{\bf Definition}
\newtheorem*{conjecture}{\bf Conjecture}
\newtheorem{example}{\bf Example}
\newtheorem*{example*}{\bf Example}
\theoremstyle{remark}
\newtheorem{remark}{\bf Remark}
\newtheorem*{remark*}{\bf Remark}
\newtheorem*{cremark}{\bf Remark on the co-orientation of the elliptic domain}
\newtheorem*{Aremark}{\bf Affine Remark}
\newtheorem*{Eremark}{\bf Euclidean Remark}
\newtheorem*{Lremark}{\bf Legendrian Remark}
\maketitle

\noindent
{\small \bf Abstract. \rm 
We show some generic (robust) properties of smooth surfaces immersed in the real $3$-space 
(Euclidean, affine or projective), in the neighbourhood of a {\em godron} (called also 
{\em cusp of Gauss}): an isolated 
parabolic point at which the (unique) asymptotic direction is tangent to the parabolic 
curve. With the help of these properties and a projective invariant that we associate 
to each godron we present all possible local configurations of the {\em flecnodal curve} 
at a generic swallowtail in $\R^3$. We present some global results, for instance: 
{\em In a hyperbolic disc of a generic smooth surface, the flecnodal curve has 
an odd number of transverse self-intersections (hence at least one self-intersection).}
\medskip

\smallskip

}


\section{Introduction}
A generic smooth surface in $\R^3$ has 
three (possibly empty) parts: 
an open {\em hyperbolic domain} at which 
the Gaussian curvature $K$ is negative,
an open {\em elliptic domain} 
at which 
$K$ is positive
and a {\em parabolic curve} at which 
$K$ vanishes. 
A {\em godron} is a parabolic point at which the (unique) asymptotic 
direction is tangent to the parabolic curve. 
We present various robust geometric properties of generic 
surfaces, associated to the godrons. For example 
(Theorem~\ref{godron_flattening}):\newline

\noindent
{\em Any smooth curve of a surface of $\R^3$ tangent to the parabolic curve 
at a godron $g$ has at least $4$-point contact with the tangent plane of the 
surface at $g$.}
\medskip

The line formed by the inflection 
points of the asymptotic curves in the hyperbolic domain is called 
{\em flecnodal curve}. The next theorem is well known. 

\begin{theorem}\label{tangency} 
{\rm (\cite{Salmon, Kortewegpp, Platonova, Landis, Banchoff})} 
At a godron of a generic smooth surface the flecnodal curve is (simply) 
tangent to the parabolic curve.
\end{theorem}

For any generic smooth surface we have the following global result 
(Proposition~\ref{euler=2} and Theorem~\ref{8}): 
\medskip

\noindent
{\em A closed parabolic curve bounding a hyperbolic disc has a positive 
even number of godrons, and the flecnodal curve lying in that disc has 
an odd number of transverse self-intersections} (thus at least one self-intersection 
point). 
\medskip 

The {\em conodal curve} of a surface $S$ is the closure 
of the locus of points of contact of $S$ with its {\em bitangent planes} 
(planes which are tangent to $S$ at least at two distinct points). 
It is well known (\cite{Salmon, Kortewegpp}) that: 
\medskip

\noindent
{\em At a godron of a generic smooth surface the conodal curve is (simply) 
tangent to the parabolic curve.}
\medskip

So the parabolic, flecnodal and conodal curves of a surface are mutually 
tangent at the godrons. At each godron, these three tangent curves determine 
a projective invariant $\rho$, as a cross-ratio (see the cr-invariant below). We 
show all possible configurations of these curves at a godron, according to the value of 
$\rho$ (Theorem~\ref{rho-classification}). There are six generic configurations, 
see Fig.~\ref{classification}. 

The invariant $\rho$ and the geometric properties of the godrons presented 
here are useful for the study 
of the local affine (projective) differential properties of swallowtails. 
So, for example, we present all generic configurations of the flecnodal curve 
in the neighbourhood of a swallowtail point of a surface of $\R^3$ in general 
position (see Theorem~\ref{4-fronts} -- Fig.~\ref{duality} and Theorem~\ref{7-fronts} -- 
Fig.~\ref{7_fronts_fig}). 
\medskip

Our results are related to several mathematical theories as, for instance, 
implicit differential equations (Davidov \cite{Davidov}), contact geometry and 
Legendrian singularities (Arnold \cite{arnoldcw, Arnoldtpwp}), and more closely to 
differential geometry, singularities of projections and tangential singularities 
(Bruce, Giblin, Tari \cite{B_G_T_95, B_G_T_98}, 
Goryunov \cite{goryunov83}, Landis \cite{Landis}, Platonova \cite{Platonova}, 
Banchoff, Thom \cite{b-Thom}).
\medskip

The paper is organised as follows. In section \ref{background}, we recall  
the classification of points of a generic smooth surface in terms of the 
order of contact of the surface with its tangent lines. In section~\ref{results}, 
we give some definitions and present our results, proving some of them directly. 
Finally, in section \ref{proofs}, we give the proofs of the theorems. 
\medskip

\noindent
{\small{\bf Acknowledgements}.\ ~
I would like to thank S.~Janeczko, W.~Domitrz and the Banach Centre for 
their hospitality and for the nice environment to do mathematics, to F.~Aicardi 
and D.~Meyer for useful comments and to E.~Ghys and D.~Serre for the references 
\cite{Kortewegpp,Korteweggtp,Levelt}.}


\section{Projective properties of smooth surfaces}\label{background}

The points of a generic smooth surface in the real 3-space (projective, 
affine or Euclidean) are classified in terms of the contact 
of the surface with its tangent lines. In this section, we recall this 
classification and some terminology. 

A generic smooth surface $S$ is divided in three (possibly empty) parts: \newline
(E) An open {\em domain of elliptic points}: 
there is no real tangent line exceeding 2-point contact with $S$; \ \newline 
(H) An open 
{\em domain of hyperbolic points}: there are 
two such lines, called {\em asymptotic lines} (their directions at the 
point of tangency are called {\em asymptotic directions}); and \ \ \newline
(P) A smooth {\em curve of parabolic points}: a unique, but double, 
asymptotic line. 

The {\em parabolic curve}, divides $S$ into the {\em elliptic} and {\em hyperbolic domains}. 

In the closure of the hyperbolic domain there is: \newline
(F) A smooth immersed {\em flecnodal curve}: it is formed by the points at which an asymptotic 
tangent line exceeds 3-point contact with $S$.  
 
\begin{figure}[ht]
\centerline{\psfig{figure=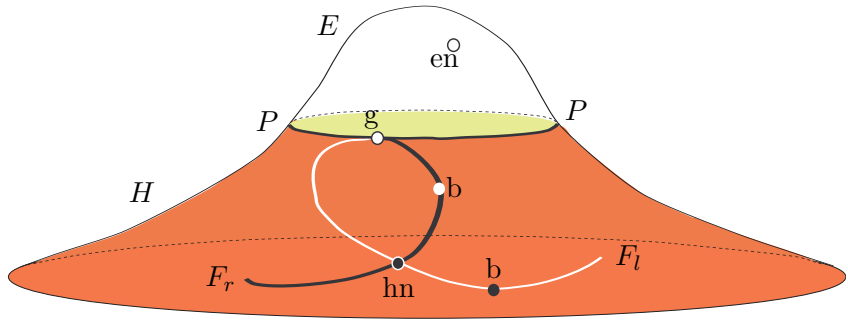,height=4.2cm}}

\begin{picture}(0,0)

\put(71,58){$H$}
\put(119,85){$P$}
\put(142,121){$E$}
\put(167,22){hn} 
\put(160,88){g}
\put(185,109){en}
\put(191,60){b}
\put(206,29){b}
\put(236,89){$P$}
\put(255,34){$F_l$} 
\put(100,26){$F_r$} 
\end{picture}

\caption{\small The 8 tangential singularities of a generic smooth surface.}
\label{points}
\end{figure}

One may also encounter isolated points of the following four types: \\
(g) A {\em godron} is a parabolic point at which the 
(unique) asymptotic direction is tangent to the parabolic curve; 
(hn) A {\em hyperbonode} is a point of the simplest self-intersection 
of the flecnodal curve; (b) A {\em biflecnode} is a point of the flecnodal curve 
at which one asymptotic tangent exceeds $4$-point contact with $S$ 
(it is also called {\em biinflection}) ;
(en) An {\em ellipnode} is a real point in the elliptic domain
of the simplest self-intersection of the complex conjugate flecnodal curves 
associated to the complex conjugate asymptotic lines. In Fig.~\ref{points} the hyperbolic domain 
is represented in gray colour and the elliptic one in white. The flecnodal curve 
has a left branch $F_l$ (white) and a right branch $F_r$ (black). These branches will be 
defined in the next section.

The term ``godron'' is due to R.~Thom \cite{k-Thom}. In other papers one can find the terms  
``special parabolic point'' or ``cusp of the Gauss map''. We keep Thom's 
terminology since it is shorter. 
Here we will study the local projective differential properties of the godrons. 

The above 8 classes of tangential singularities, Theorem~\ref{tangency} and all the theorems 
presented in this paper are projectively invariant and are robust features of a smooth 
surface, that is, they are stable in the sense that under a sufficiently small 
perturbation (taking derivatives into account) they do not vanish 
but only deform slightly. Seven of these classes were known at the end of the 19th 
century in the context of the enumerative geometry of complex algebraic surfaces, 
with prominent works of Cayley, Zeuthen and Salmon, see \cite{Salmon}. 
For these seven classes, the normal forms of surfaces at such points up to the 5-jet, 
under the group of projective transformations, 
were independently found by E.E. Landis (\cite{Landis}) and 
O.A. Platonova (\cite{Platonova}). The ellipnodes 
were found by D.~Panov (\cite{Dima}) who called them  special elliptic points.

For surfaces in $\R^3$, these tangential singularities depend 
only on the affine structure of $\R^3$ (because they depend only on the contact 
with lines), that is, they are independent of any 
Euclidean structure defined on $\R^3$ and of the Gaussian curvature of the 
surface which could be induced by such a Euclidean structure. 

Another definition of godron (cusp of Gauss) is in terms of the contact 
with the tangent plane. In this setting, a useful tool for analysis is the so 
called `height function' (cf. \cite{B_G_T_95, B_G_T_98}), which can be defined once 
some Euclidean structure is fixed. Anyway, the singularities of the contact with the 
tangent plane `are expressed geometrically', independently of any Euclidean structure,  
by the singularities of the dual surface: 
An ordinary godron corresponds to a swallowtail point of the dual surface, 
that is, to an $A_3$ (Legendre) singularity (cf. \cite{avg, Arnoldtpwp}). There are 
also double (unstable) godrons, corresponding to an $A_4$ bifurcation 
(two swallowtails born or dying), that is, to an $A_4$ (Legendre) singularity 
(see \S\ref{degenerated_rho=1}). 

We will say that a godron is {\em simple} if it corresponds to a swallowtail point 
of the dual surface. All godrons of a surface in general position are simple. 

Both types of tangential 
singularities (contact with lines and contact with planes) were extensively 
used by Cayley, Zeuthen and Salmon, see \cite{Salmon}. 


Besides the smooth surfaces, we also consider surfaces admitting wave front 
singularities (section \ref{fronts}) and 
we study the behaviour of the flecnodal curve near the swallowtail points. 


\section{Statement of results}\label{results}
Consider the pair of fields of asymptotic directions in the hyperbolic domain. An 
{\em asymptotic curve} is an integral curve of a field of asymptotic 
directions.  
\medskip

\noindent
{\bf Left and right asymptotic and flecnodal curves}.\ ~
Fix an orientation in the $3$-space $\RP^3$ (or in $\R^3$).
The two asymptotic curves passing through a point of the hyperbolic 
domain of a generic smooth surface can be distinguished in a natural 
geometric way: One twists like a left screw and the other like a right screw. 
More precisely, a regularly parametrised smooth curve is said to be a 
{\em left (right) curve} if its first three derivatives at each point form 
a negative (resp. a positive) frame. 

\begin{proposition}\label{left-right}
At a hyperbolic point of a surface one asymptotic curve is left and the 
other one is right.
\end{proposition}

A proof is given (for generic surfaces) in Euclidean Remark below. 

The hyperbolic domain is therefore foliated by a family of left asymptotic curves 
and by a family of right asymptotic curves. The corresponding asymptotic tangent 
lines are called respectively {\em left} and {\em right asymptotic lines}.

\begin{definition*}
The {\em left} ({\em right}) {\em flecnodal curve} $F_l$ (resp. $F_r$) of a surface 
$S$ consists of the points of the flecnodal curve of $S$ whose asymptotic line, 
having higher order of contact with $S$, is a left (resp. right) asymptotic line.  
\end{definition*}

The following statement (complement to Theorem~\ref{tangency}) is used and implicitly 
proved (almost explicitly) in \cite{Uribetesis, surf-evolution}. A proof is given in 
section~\ref{proofs}, see Fig.~\ref{points}: 

\begin{theorem}\label{separates}
A simple godron separates locally the flecnodal curve into its right and left branches.
\end{theorem}

\begin{definition*}
A {\em flattening} of a generic curve is a point at which the first three derivatives 
are linearly dependent. Equivalently, a flattening is a point at which the curve 
has at least $4$-point contact with its osculating plane. 
\end{definition*} 

The flattenings of a generic curve are isolated points separating the 
right and left intervals of that curve. 

\begin{Eremark}
If we fix an arbitrary Euclidean structure in the affine oriented space 
$\R^3$, then the lengths of the vectors and the angles between vectors are 
defined. Therefore, for such Euclidean structure, the torsion $\tau$ of a 
curve and the Gaussian curvature $K$ of a surface are defined. In this case 
{\em a point of a curve is 
right, left or flattening if the torsion at that point satisfies $\tau>0$, $\tau<0$ or $\tau=0$, 
respectively}. The Gaussian curvature $K$ on the hyperbolic domain of a smooth surface 
is negative. The Beltrami-Enepper Theorem states that 
{\em the values of the torsion of the two asymptotic curves passing through a 
hyperbolic point with Gaussian curvature $K$ are given by $\tau = \pm \sqrt {-K}$}. 
This proves Proposition~\ref{left-right}. 
\end{Eremark}

\begin{definition*}
An {\em inflection} of a (regularly parametrised) smooth curve is a point at which 
the first two derivatives are linearly dependent. Equivalently, an inflection is a 
point at which the curve has at least $3$-point contact with its tangent line. 
\end{definition*} 

A generic curve in the affine space $\R^3$ has no inflection. However, a generic 
$1$-parameter family of curves can have isolated parameter values for which the 
corresponding curve has one isolated inflection. 

\begin{theorem}\label{godron_flattening}
Let $S$ be a smooth surface. All smooth curves of $S$ which are 
tangent to the parabolic curve at a godron $g$ have either a flattening or an 
inflection at $g$, and their osculating plane is the tangent
plane of $S$ at $g$.  
\end{theorem}

The proof of Theorem~\ref{godron_flattening} is given in section \ref{proofs}.

\begin{corollary}
If a point is a godron of a generic smooth surface, then it is a flattening of 
both the parabolic curve and the flecnodal curve. 
\end{corollary}

\begin{remark*}
The converse is not true: A flattening of the parabolic curve or of the flecnodal curve 
is not necessarily a godron.  
\end{remark*}

\subsection{The cr-invariant and classification of godrons}\label{sec:cr-invariant}

{\bf The conodal curve}.\ ~ 
Let $S$ be a smooth surface. A {\em bitangent plane} of $S$ 
is a plane which is tangent to $S$ at least at two distinct points 
(which form a pair of {\em conodal points}). 
The {\em conodal curve} $D$ of a surface $S$ is the closure of the 
locus of points of contact of $S$ with its bitangent planes. 
\medskip

At a godron of $S$, the curve $D$ is simply tangent 
to the curves $P$ (parabolic) and $F$ (flecnodal). This fact will be clear 
from our calculation of $D$ for Platonova's normal form of godrons. 
\medskip

\noindent
{\bf The projective invariant}.\ ~ 
At any simple godron $g$, there are three tangent smooth curves $F$, $P$ and $D$, 
to which we will associate a projective invariant: \\

Consider the Legendrian 
curves $L_F$, $L_P$, $L_D$ and $L_g$ (of the $3$-manifold of contact 
elements of $S$, $PT^*S$) consisting of the contact elements of $S$ tangent 
to $F$, $P$, $D$ and to the point $g$, respectively (the contact elements of $S$ tangent 
to a point are just the contact elements of $S$ at that point, that is, $L_g$ 
is the fibre over $g$ of the natural projection $PT^*S \rightarrow S$). These four Legendrian 
curves are tangent to the same contact plane $\Pi$ of $PT^*S$. The tangent 
directions of these curves determine four lines $\ell_F$, $\ell_P$, $\ell_D$ and $\ell_g$, 
through the origin of $\Pi$. 

\begin{definition*}\label{invariant}
The {\em cr-invariant} $\rho(g)$ of a godron $g$ is defined as the 
cross-ratio of the lines $\ell_F$, $\ell_P$, $\ell_D$ and $\ell_g$ of $\Pi$: 
$$\rho(g)=(\ell_F,\ell_P,\ell_D,\ell_g).$$ 
\end{definition*} 

\noindent
{\bf Platonova's normal form}.\ ~
According to Platonova's Theorem \cite{Platonova}, in the neighbourhood 
of a godron, a surface can be sent by projective transformations to the 
normal form
$$z=\frac{y^2}{2}-x^2y+\lambda x^4+\f(x,y) \ \ \ \mbox{(for some $\lambda\neq 0,\frac{1}{2}$)} \eqno(G1)$$
where $\f$ is the sum of homogeneous polynomials in $x$ and $y$ of degree 
greater than $4$ and (possibly) of flat functions. 

\begin{theorem}\label{rho=2l}
Let $g$ be a godron, with cr-invariant value $\rho$, of a generic smooth surface $S$. 
Put $S$ (after projective transformations) in Platonova's normal form $(G1)$. 
Then the coefficient $\lambda$ equals $\rho/2$.
\end{theorem}

It turns out that among the $2$-jets of the curves in $S$, tangent to $P$ at 
a godron, there is a special $2$-jet at which ``something happens''. 
We introduce it in the following lemma. 
\medskip

\noindent
{\bf Tangential Map and Separating $2$-jet}.\ ~ Let $g$ be a godron of a 
generic smooth surface $S$. The {\em tangential map} of $S$, 
$\tau_S:S\rightarrow (\RP^3)^\vee$, associates to each point of $S$ 
its tangent plane at that point. The image $S^\vee$ of $\tau_S$ is called 
the {\em dual surface of $S$}. 

Write $J^2(g)$ for the set of all $2$-jets of curves of $S$ tangent to 
$P$ at $g$. By {\em the image of a $2$-jet $\g$ in $J^2(g)$ under the 
tangential map $\tau_S$} we mean the image, under $\tau_S$, of any curve 
of $S$ whose $2$-jet is $\g$. By Theorem \ref{godron_flattening}, 
all the $2$-jets of $J^2(g)$ (and also the $3$-jets of curves on $S$ 
tangent to $P$ at $g$) are curves lying in the tangent plane of $S$ 
at $g$. In suitable affine coordinates, 
the elements of $J^2(g)$ can be identified with the curves 
$t\mapsto (t,ct^2,0)$, $c\in \R$. 

\begin{slemma}\label{separating}
There exists a unique $2$-jet $\s$ in $J^2(g)$ {\rm (that we call 
{\em separating $2$-jet at $g$})} satisfying the following properties: 
\medskip

\noindent
$(a)$ The images, under $\tau_S$, of all elements of $J^2(g)$ different from 
$\s$ are cusps of $S^\vee$ sharing the same tangent line $\ell_g^\vee$, at $\tau_S(g)$. 
\medskip

\noindent
$(b)$ The image of $\s$ under $\tau_S$ is a singular curve of $S^\vee$ whose 
tangent line at $\tau_S(g)$ is different from $\ell_g^\vee$. 
\medskip

\noindent
$(c)$ {\rm (separating property)}: The images under $\tau_S$ of any two elements of 
$J^2(g)$, separated by $\s$, are cusps pointing in opposite directions.  
\end{slemma}


\begin{remark*}\label{y=x^2}
Once a godron with cr-invariant $\rho$ of a smooth surface is sent 
(by projective transformations) to the normal form 
$z=y^2/2-x^2y+\rho x^4/2+\f(x,y)$, the separating $2$-jet is independent of $\rho$:  
It is given by the equation $y=x^2$, in the $(x,y)$-plane. 
\end{remark*}




For almost all values of 
$\rho$ the curves $F$, $P$ and $D$ are simply tangent one to the others. 
However, for isolated values of $\rho$ two of these curves may have higher 
order of tangency and then some bifurcation occurs. We will look for the values 
of $\rho$ at which `something happens'. 
\medskip

\noindent
{\bf Canonical coefficients}.\ 
Consider a godron $g$ with cr-invariant $\rho$ and suppose that the surface was 
sent (by projective transformations) to the normal form 
$z=y^2/2-x^2y+\rho x^4/2+\f(x,y)$. 
The plane curves $\bar F$, $\bar P$ and $\bar D$, which are the projections 
of $F$, $P$ and $D$ to the $(x,y)$-plane along the $z$-axis, have the same $2$-jet as 
$F$, $P$ and $D$, respectively (since, according to Theorem~\ref{godron_flattening}, $F$, $P$ and $D$ 
have at least $4$-point contact with the $(x,y)$-plane).
These $2$-jets correspond to three parabolas $~y=c_\sF x^2$, $~y=c_\sP x^2~$ and 
$~y=c_\sD x^2,~$ whose coefficients $~c_\sF$, $~c_\sP~$ and $~c_\sD~$ we  
call the {\em canonical coefficients} of the curves $F$, $P$ and $D$, respectively.  

\begin{figure}[ht]
\centerline{\psfig{figure=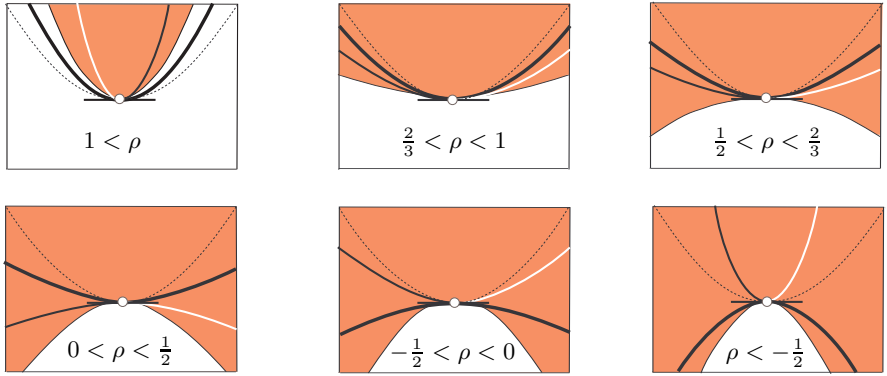,height=5cm}}

\begin{picture}(0,0)

\put(47,100){\small\bf $1<\rho$} 
\put(167,100){\small\bf $\frac{2}{3}<\rho<1$} 
\put(285,100){\small\bf $\frac{1}{2}<\rho<\frac{2}{3}$} 
\put(41,21){\small\bf $0<\rho <\frac{1}{2}$} 
\put(163,20){\small\bf $-\frac{1}{2}<\rho<0$} 
\put(290,20){\small\bf $\rho<-\frac{1}{2}$} 
\end{picture}

\caption{\small The configurations of the curves $F$ (half-white half-black curves), 
$P$ (boundary between white and gray domains), $D$ (thick curves), the separating 
$2$-jet (broken curves) and the asymptotic line (horizontal segments) at generic godrons.}
\label{classification}
\end{figure}

The configuration of the curves $F$, $P$ and $D$ 
with respect to the asymptotic line and the separating $2$-jet at $g$ 
is equivalent to the configuration at the origin of the curves 
$\bar F$, $\bar P$ and $\bar D$ with respect to the parabolas $y=0\cdot x^2=0$ and $y=1\cdot x^2$ 
on the $(x,y)$-plane (see the above Remark).  This configuration is determined by the 
relative positions of the canonical coefficients $c_\sF$, $c_\sP$, $c_\sD$, with respect to 
the numbers $c_\s=1$ and $c_{al}=0$, in the real line:

\begin{theorem}\label{rho-classification}
Given a simple godron $g$ of a smooth surface,   
there are six possible configurations of the curves $F$, $P$ and $D$ 
with respect to the separating $2$-jet and to the asymptotic line at $g$ 
{\rm (they are represented in Fig.~\ref{classification})}. The actual configuration 
at $g$ depends on which of the following six open intervals 
the cr-invariant $\rho(g)$ belongs to, respectively: 
$$
\begin{array}{ccc}
\rho\in(1,\infty) &  \iff &  1<c_\sD<c_\sP<c_\sF;\\ 
\rho\in(\frac{2}{3},1) &  \iff &  0<c_\sP<c_\sF<c_\sD<1;\\
\rho\in(\frac{1}{2},\frac{2}{3}) &  \iff &  c_\sP<0<c_\sF<c_\sD<1;\\
\rho\in(0,\frac{1}{2})  & \iff &  c_\sP<c_\sF<0<c_\sD<1;\\
\rho\in(-\frac{1}{2},0) &  \iff &  c_\sP<c_\sD<0<c_\sF<1;\\
\rho\in(-\infty,-\frac{1}{2}) &  \iff &  c_\sP<c_\sD<0<1<c_\sF.
\end{array}
$$
\end{theorem}

Besides the 5 exceptional values of $\rho$, given in Theorem~\ref{rho-classification}, 
we will present separately (\S\ref{g_contour}) other important exceptional values of $\rho$.

\subsection{The index of a godron}\label{index-definition}

\begin{figure}[ht]
\begin{picture}(0,0)
\put(123,5){\small\bf $(+)$} 
\put(237,5){\small\bf $(-)$} 
\end{picture}
\centerline{\psfig{figure=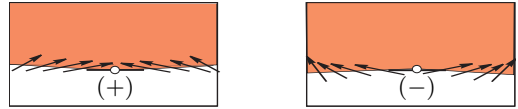,height=1.4cm}}
\caption{\small A positive godron and a negative godron.}
\label{index}
\end{figure}

\begin{definition*}
A godron is said to be {\em positive} or {\em of index $+1$} 
(resp.~{\em negative} or {\em of index $-1$}) if at the neighbouring parabolic 
points the half-asymptotic lines, directed to the hyperbolic domain, point 
towards (resp.~away from) the godron \ -- Fig.~\ref{index} (some authors use the term 
{\em hyperbolic} (resp. {\em elliptic})).   
\end{definition*}

\noindent
{\bf The asymptotic double of the hyperbolic domain}.\ ~ A godron $g$ can be positive or 
negative, depending on the index of the direction field, which is naturally 
associated to $g$, on the {\em asymptotic double} $\A$ of $S$: 
The {\em asymptotic double} of $S$ is the surface $\A$ in the manifold of contact elements 
of $S$, $PT^*S$, consisting of the field of asymptotic directions. 
It doubly covers the hyperbolic domain, and its projection to $S$ has a 
fold singularity over the parabolic curve. There is an
{\em asymptotic lifted field of directions} on the surface $\A$, constructed in the following way. 
At each point of the contact manifold $PT^*S$ a contact plane is applied, 
in particular at each point of $\A$.  
Consider a point of the smooth surface $\A$ and assume that the tangent 
plane of $\A$ at this point  does not coincide with the contact 
plane. Then these two planes intersect along a straight line tangent to 
$\A$. The same holds at all nearby points in $\A$. This defines a 
smooth direction field on $\A$ which vanishes only at the points 
where those planes coincide: over the godrons.

\begin{figure}[ht]
\centerline{\psfig{figure=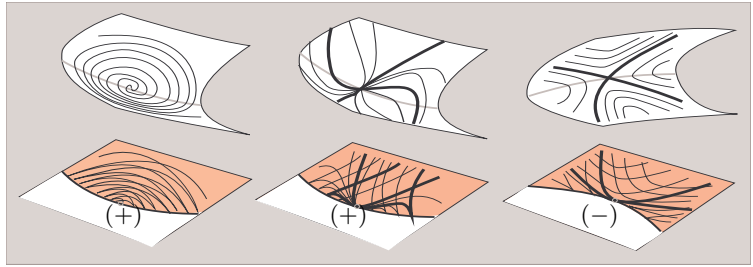,height=3.5cm}}

\begin{picture}(0,0)
\put(80,30){\small\bf $(+)$} 
\put(165,30){\small\bf $(+)$} 
\put(260,30){\small\bf $(-)$} 
\end{picture}

\caption{\small The asymptotic double of the hyperbolic domain near a godron.}
\label{double}
\end{figure}

If $g$ is a positive godron, then the index of this direction field at 
its singular point equals $+1$, the point being a node or a focus; if 
$g$ is negative, the index equals $-1$ and the point is a saddle. See Fig.~\ref{double}.

\begin{proposition}\label{positive-rho=1}
A godron $g$ is positive (negative) if and only if the value of its cr-invariant 
$\rho$ satisfies: $\rho(g)>1$ (resp. $\rho(g)<1$). 
\end{proposition}

We say that the elliptic domain of a smooth surface $S$ 
is {\em locally convex in the neighbourhood of a godron} if, when projected to 
the tangent plane to $S$ (from any point, exterior to this plane), the image 
of the elliptic domain is locally convex: The tangent line of the parabolic 
curve being locally inside the image of the hyperbolic domain (the hyperbolic 
domain being locally convex if this line lies locally in inside the image of the 
elliptic domain). 

\begin{corollary}[of Theorem~\ref{rho-classification}]\label{2/3}
At a godron $g$ with cr-invariant $\rho$ the hyperbolic (elliptic) domain is 
locally convex if and only if $\rho >2/3$ (resp. $\rho <2/3$).
\end{corollary}

\begin{proof}
The hyperbolic (elliptic) domain is locally convex at $g$ if and only if the canonical 
coefficient $c_\sP$ of the parabolic curve is positive (resp. negative). 
\end{proof}

\begin{theorem}\label{franca} \

\noindent
$(a)$ In the neighbourhood of any positive godron the hyperbolic domain is locally 
convex. 
\medskip

\noindent
$(b)$ There exist negative godrons for which the neighbouring hyperbolic 
domain is locally convex.
\medskip

\noindent
$(c)$ At the negative godrons for which the neighbouring hyperbolic 
domain is locally convex, the flecnodal curve lies locally between $P$ and $D$
{\rm (see Fig.~\ref{classification})}. Moreover, the cr-invariant satisfies: 
$\frac{2}{3}<\rho <1$. 
\end{theorem}

\begin{proof}
The theorem follows from Proposition \ref{positive-rho=1} and Corollary~\ref{2/3}. 
\end{proof}

Items $(a)$ and $(b)$ of Theorem~\ref{franca} were discovered by F.~Aicardi \cite{Franca}. 
\bigskip

\begin{corollary}\label{cubic}
All godrons of a cubic surface in $\RP^3$ are negative.
\end{corollary}

\begin{proof}
By the definitions of asymptotic curve and of flecnodal curve, any straight line contained 
in a smooth surface is both an asymptotic curve and a connected component of the 
flecnodal curve of that surface. 

Let $S$ be an algebraic surface of degree $3$. At a point of the flecnodal 
curve, an asymptotic line has at least 
$4$-point contact with $S$. Since $S$ is a cubic surface, this line must lie completely 
in $S$. So the flecnodal curve of $S$ consists of straight lines. 

At a godron $g$ of $S$, the tangent line to the parabolic 
curve (that is, the flecnodal curve) lies in the hyperbolic domain. Thus the 
neighbouring elliptic domain is locally convex. Therefore, by the above theorem, 
$g$ is negative. 
\end{proof}

\noindent
{\bf Factorisable polynomials}. \ 
A set $\{\ell_1, \ldots,\ell_n\}$ of real affine functions on the plane
is said to be {\em in general position} if: (i) The lines $\ell_i=0$, $\ell_j=0$ 
are not parallel ($i\neq j$, $i,j\in\{1,\ldots,n\}$), and (ii) For any 
$i\in\{1,\ldots,n\}$,  the line $\ell_i=0$ contains no critical point of 
the function $\prod_{j\neq i}\ell_j$. 
The product $\prod \ell_i$ of $n$ real affine functions in general position 
is called a {\em factorisable polynomial}. 

In \cite{adriana}, A. Ortiz-Rodr\'{\i}guez 
proved, among other things, that for any real factorisable polynomial of degree $n$, 
$f=\prod \ell_i$, the following holds: (i) {\em The lines $\ell_i$ are the only 
components of the flecnodal curve 
of the graph of $f$}, (ii) {\em This graph has exactly $n(n-2)$ godrons}, and  

\begin{proposition}\label{all_negative}
{\rm (Theorem 1 and Lemma 13 of \cite{adriana})}
All godrons of the graph of a real factorisable polynomial are negative.
\end{proposition}

\begin{proof}
Theorem~\ref{franca} provides an alternative and very simple proof 
of Proposition~\ref{all_negative}: Since the flecnodal curve consists of 
straight lines, at each godron $g$ the asymptotic tangent line and the 
asymptotic curve coincide with one of such straight lines. 
Thus, the elliptic domain is locally convex at $g$, and hence $g$ is a 
negative godron (by Theorem~\ref{franca}).
\end{proof} 

\subsection{Locating the left and right branches of the flecnodal curve}

\begin{cremark}
Each connected component of the elliptic domain is `naturally' co-oriented: 
At each elliptic point the surface lies locally on one of the two half-spaces 
determined by its tangent plane at that point. This half-space, that we 
name {\em positive half-space}, determines a {\em natural co-orientation} 
on each connected component of the elliptic domain. By continuity, the natural 
co-orientation extends to the parabolic points (where the parabolic curve 
is smooth). At the parabolic points a positive 
half-space is therefore also defined. 
\end{cremark}

This observation has strong topological consequences. For example: 

\begin{theorem}
The elliptic domain of any smooth surface in the $3$-space (Euclidean, affine or 
projective) can not contain a M\"obius strip. 
\end{theorem}

\begin{proof}
If a M\"obius strip $M$ were contained in the elliptic domain $E$ of a surface, then 
it would be contained in a connected component of $E$, since $M$ is connected. Now, 
the theorem follows since each connected component of $E$ has a natural co-orientation 
(and $M$ is not co-orientable). 
\end{proof}

In the neighbourhood of a godron $g$ of a smooth surface $S$, we can distinguish 
explicitly which branch of the flecnodal curve is the right branch and which 
is the left one. For this, we need only to know the index of $g$ and the 
natural co-orientation of $S$ (given by the positive half-space at $g$):

Let $g$ be a godron, with $\rho\neq 1$, of a smooth surface $S$. 
Take an affine coordinate system $x,y,z$ such that the $(x,y)$-plane is 
tangent to $S$ at $g$, and the $x$-axis is tangent to the parabolic curve 
at $g$ (thus also tangent to $F$ at $g$). Direct the positive $z$-axis 
to the positive half-space at $g$. Direct the positive $y$-axis 
towards the neighbouring hyperbolic domain. Finally, direct the positive 
$x$-axis in such way that any basis $(e_x,e_y,e_z)$ of 
$x,y,z$ form a positive frame for the fixed orientation 
of $\R^3$ (or of $\RP^3$). 

So one can locally parametrise the flecnodal curve at $g$ by projecting it to the $x$-axis. 

\begin{theorem}\label{lr}
Under the above parametrisation, 
the left and right branches of the flecnodal curve at $g$ correspond locally to the 
negative and positive semi-axes of the $x$-axis, respectively, if and only if $g$ 
is a positive godron. The opposite correspondence holds for a negative godron. 
\end{theorem}

In other words, if you stand on the tangent plane of $S$ at $g$ in the positive 
half-space and you are looking from the elliptic domain to the hyperbolic one,
then you see the right (left) branch of the flecnodal curve on your right hand 
side if and only if $g$ is a positive (resp. negative) godron. So the index of $g$ 
determines and is determined by the side on which the right branch of $F$ is 
located.

\begin{remark*}
Theorems~\ref{separates} and \ref{lr} (which are local theorems)
together with the natural co-orientation of the elliptic domain, are the key 
elements to prove the global theorem (Theorem~\ref{8}) of section \ref{discs}. 
They imply that some (global) configurations of the flecnodal 
curve are forbidden. So, for example, there is no surface having a hyperbolic 
disc without hyperbonodes. 
\end{remark*}

\subsection{The flec-godrons: Degenerated godrons with $\rho=0$}\label{degenerated_rho=0}

After the preceding sections, a natural question arises: 
What happens if the cr-invariant equals $0$ or $1$? 

The godrons for which the cr-invariant equals $0$ or $1$ are 
degenerated godrons. 
We will explain the meaning of these degeneracies and 
describe the behaviour of such degenerated godrons under a small 
perturbation of the surface inside a generic one parameter family 
of smooth surfaces. 
\medskip

\noindent
{\bf The case $\rho=0$}. \ 
If $\rho=0$ then the $4$-jet given by the (Platonova) normal form, used above, 
defines just a cubic surface, which is absolutely not generic: The asymptotic 
tangent line at that godron has `infinite point-contact' with the surface and 
coincides with the flecnodal curve. In order to understand the behaviour 
of the flecnodal curve (and  the geometric properties of the surface)  
at the godrons with $\rho=0$, we need to add some terms of degree $5$, which in this 
case are relevant and break the symmetry. In fact, we have the 

\begin{proposition}
At a godron of a smooth surface the asymptotic tangent line has 
$4$-point contact with the surface if and only if $\rho\neq 0$.  
At a godron with $\rho=0$ the asymptotic tangent line and the surface have at 
least $5$-point contact. 
\end{proposition}

\begin{proof}
Consider the godron $g$ with cr-invariant $\rho$ (at the origin) of the surface $S$ 
given by $z=\frac{y^2}{2}-x^2y+ \frac{\rho}{2}x^4 + \f(x,y)$, 
where $\f(x,y)$ is the sum of homogeneous polynomials in 
$x$ and $y$ of degree greater than $4$. Since the asymptotic tangent line $\ell$ 
at $g$ is the $x$-axis, its order of contact with $S$ at $g$ is the multiplicity 
of the zero of the function $g\comp\g$ at $t=0$, where 
$g(x,y,z)=-z+\frac{y^2}{2}-x^2y+ \frac{\rho}{2}x^4 + \f(x,y)$ and 
$\g(t)=(x(t),y(t),z(t))=(t,0,0)$.  

Since $g\comp \g$ has the form $(g\comp\g)(t)=\frac{\rho}{2}t^4 + at^5 +\ldots$, 
(where $a$ is the coefficient of $x^5$ in $\f(x,y)$), 
the asymptotic line $\ell$ has $4$-point contact with $S$ at $g$ if and only if 
$\rho\neq 0$ \ -- and $5$-point contact if and only if $\rho=0$ (and $a\neq 0$). 
\end{proof}

\begin{definition*}
A {\em flec-godron} of a smooth surface is a godron at which the asymptotic tangent line 
and the surface have at least $5$-point contact (being exactly $5$-point contact 
for a {\em simple flec-godron}). 
\end{definition*}

It follows that a smooth surface in general position has no flec-godron:
Under any small generic deformation of the surface the flec-godron condition 
of $5$-point contact of the surface with the asymptotic line at $g$
(or, equivalently, $\rho=0$) is destroyed. 
However, a generic $1$-parameter family of smooth surfaces can have, 
at isolated parameter values, a surface having one simple flec-godron 
The smooth surfaces 
in the $3$-space having a flec-godron form a {\em discriminant hypersurface} in the space 
of smooth surfaces. We shall describe the local bifurcation of the surface (of its tangential 
singularities) when a generic $1$-parameter family traverses this discriminant hypersurface. 

\begin{example*}
The godron of the surface $z=\frac{y^2}{2}-x^2y+ ax^5$, with $a\neq 0$, is a simple flec-godron. 
We will perturb this surface inside a generic one parameter family 
of smooth surfaces in which the parameter is the cr-invariant $\rho$: 
$$z=\frac{y^2}{2}-x^2y+ \frac{\rho}{2}x^4 + ax^5.$$

\begin{figure}[ht]
\begin{picture}(0,0)
\put(63,-10){\small $\rho<0$}
\put(174,-10){\small $\rho=0$} 
\put(286,-10){\small $\rho>0$}
\end{picture}
\centerline{\psfig{figure=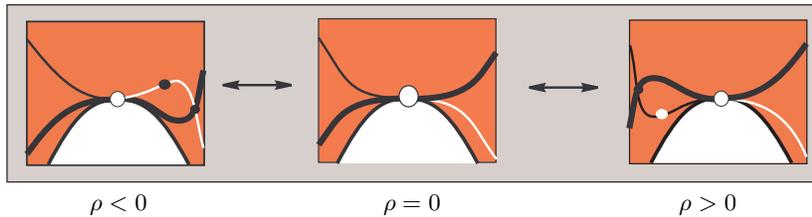,height=2.4cm}}
\caption{\small The transition of the flecnodal curve (half-black half-white curve) and of the conodal 
curve (thick curve) at a flec-godron transition: $\rho=0$.}
\label{flexgodron}
\end{figure}

The flecnodal, conodal and parabolic curves of this surface (with $a>0$) are depicted 
in Fig.~\ref{flexgodron} for $\rho<0$, $\rho=0$ and $\rho>0$ 
(compare with Fig.~\ref{classification}).  

Indeed, the power series expansions of the curves $\bar F$ and $\bar D$ start by 
$$y=\rho(2\rho-1)x^2+10a(2\rho-1)x^3+\ldots\ \ \mbox{ and }\ \   y=\rho x^2+2ax^3+\ldots,$$ 
respectively. So, for $\rho=0$ we have, $\bar F$: $y=-10ax^3+\ldots$ and 
$\bar D$: $y=2ax^3+\ldots$ (whence the bifurcation of Fig.~\ref{flexgodron}). 
\end{example*}
\smallskip

Let $S_t$, $t\in\R$, a generic one parameter family of smooth surfaces 
such that for $t=0$ the surface $S_0$ has a flec-godron $g_0$
(for instance the above family). 

As $t$ is increasing and passing through $0$, a biflecnode $b_t$ of $S_t$ 
is moving along the flecnodal curve, passing, at $g_0$, from one branch to the other of 
the flecnodal curve. Roughly speaking, `a flec-godron is the superposition of a biflecnode and 
a godron': 

\begin{theorem}\label{flec-godron_theorem}
At the flec-godron (transition) the following happens: For any $t\neq 0$, 
sufficiently close to $0$, the surface $S_t$ possesses an ordinary godron and 
a neighbouring biflecnode . The biflecnode being left for all (small) $t$ of a given sign 
and being right for all (small) $t$ of the opposite sign. Moreover, for $t=0$ both the 
flecnodal curve and the conodal curve of $S_0$ have an inflection at $g_0$. 
\end{theorem}

Note that the index of the flec-godron of $S_0$ and of the godron 
of $S_t$, for $|t|$ sufficiently small, is $-1$ (since $\rho<1$).  
\medskip

The corresponding bifurcation of the tangential singularities on the dual surface 
is described in \S\S\ref{flec-swallowtail-section} \ -- Fig.~\ref{flex-swallowtail}.

\subsection{The bigodrons: Degenerated godrons with $\rho=1$}\label{degenerated_rho=1}

\begin{figure}[ht]
\begin{picture}(0,0) 
\put(103,-10){\small $\e<0$}
\put(174,-10){\small $\e=0$} 
\put(248,-10){\small $\e>0$}
\end{picture}
\centerline{\psfig{figure=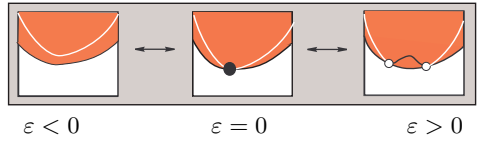,height=1.4cm}}
\caption{\small The transition at bigodron: $\rho=1$.}
\label{bigodron}
\end{figure}

\noindent
{\bf The case $\rho=1$}. \ 
If $\rho=1$, then we also have a degenerate godron, which we name {\em bigodron}: 
it is the collapse (or the birth) of two godrons with opposite indices (it is not 
a simple godron). When $\rho=1$ 
the normal form that we used above is not convenient since it is degenerate:
 $z=\frac{1}{2}(y-x^2)^2$. 
For this reason the parabolic and flecnodal curves coincide with the curve $y=x^2$ in 
the $(x,y)$-plane (this curve is sent to a point under the tangential map of $S$). 
In order to have a generic polynomial of degree four, one must add another term 
of degree four:  $z=\frac{1}{2}(y-x^2)^2\pm x^3y$. Now, the bigodron obtained is generic 
(among the bigodrons: $\rho=1$): the parabolic and flecnodal curves have $4$-point 
contact and the whole flecnodal curve is either left or right, according to the sign 
$+$ or $-$ of the term $\pm x^3y$, respectively 
(see the central part of Fig.~\ref{bigodron}). 
To understand better the geometry 
of a bigodron, we will perturb this surface inside a generic one parameter family 
of smooth surfaces: 
$$z=\frac{1}{2}(y-x^2)^2\pm x^3y+\e x^3.$$

The flecnodal and parabolic curves of this surface are depicted 
in Figure~\ref{bigodron} for $\e<0$, $\e=0$ and $\e>0$.  
When the parameter $\e$ is negative the flecnodal curve is left and does not touch 
the parabolic curve, while when $\e$ is positive the flecnodal curve touches the parabolic curve 
at two neighbouring godrons with opposite indices, and a small segment of the right flecnodal has 
appeared between these godrons. 

So, in the oriented $3$-space (Euclidean, affine or projective), there are two types of 
bigodrons: a bigodron is said to be {\em left} ({\em right}) if it 
corresponds to a bifurcation in which a small segment of the left (resp. right) branch of the 
flecnodal curve is born or vanishes. 

\begin{remark*}
The case $\rho=1$ (a bigodron) corresponds to an $A_4$ contact with the tangent plane 
(cf. \cite{B_G_T_95, B_G_T_98, arnoldcw}). That is, for the dual surface 
it corresponds to the $A_4$ bifurcation of wave fronts (two swallowtails 
are born or dying) occurring in generic one parameter families of 
fronts \cite{Arnoldwfeeml}. Thus, in the oriented $3$-space $\RP^3$ (or $\R^3$), 
there are two types of $A_4$ singularities of wave fronts. 
\end{remark*} 

\begin{remark*}
Strictly speaking,  the cr-invariant is not defined at the bigodron (or $A_4$) 
singularity, since there is no conodal curve. However, the limit value of the 
cr-invariant for both dying (or born) godrons, at the bigodron bifurcation 
moment, is equal to $1$. 
\end{remark*} 

\subsection{Elliptic discs and hyperbolic discs of smooth surfaces}\label{discs}

The following global theorem holds for any generic smooth surface:

\begin{theorem}\label{8}
In any hyperbolic disc bounded by a Jordan parabolic curve, there is an 
odd number of hyperbonodes (hence at least one). 
\end{theorem}

\begin{remark*}
The hyperbonodes and the ellipnodes of a smooth surface are the points at which 
that surface is better approximated by a quadric. Indeed, a surface in $\RP^3$, 
different from a plane, is a one sheet hyperboloid (or an ellipsoid) if and only 
if all its points are hyperbonodes (resp. ellipnodes). These points play an essential 
rôle in (and are necessary for) several bifurcations of the parabolic curve and of wave 
fronts \cite{surf-evolution}. 
\end{remark*}

The cubic surfaces in $\RP^3$ provide examples of surfaces having elliptic discs 
whose bounding parabolic curves have $0$, $1$, $2$ or $3$ negative godrons: According 
to Segre \cite{Segre}, a generic cubic surface diffeomorphic to the projective plane 
contains four parabolic curves (each one bounding an elliptic disc) and six godrons. 
According to \cite{Arnoldtpwp}, Shustin has proved that the distribution of the godrons among the 
four parabolic curves is $6=0+1+2+3$. By Corollary~\ref{cubic}, all these godrons are 
negative. 

There exist smooth surfaces having an elliptic disc whose bounding parabolic 
curve has $4$ negative godrons: 

\begin{example*}
The algebraic surface given by the equation 
$$z=(x^2-1)(y^2-1)$$ 
has an elliptic disc whose bounding parabolic curve contains $4$ godrons, all 
negatives. 
\end{example*} 

\noindent
{\bf Problem}. \  Exist there smooth surfaces in the 3-space (affine or projective) 
having an elliptic disc whose bounding parabolic curve has more than $4$ godrons, 
all of them negative? 

For a parabolic curve bounding a hyperbolic disc the situation is more restrictive: 

\begin{proposition}\label{euler=2}
The sum of the indices of the godrons on the parabolic curve bounding a hyperbolic 
disc (of a generic surface) equals two. In particular, such parabolic curve contains 
a positive even number of godrons. 
\end{proposition}

\begin{proof}
Write $H$ for the closure of the hyperbolic disc. The asymptotic double $\A$ 
is a sphere. Its Euler characteristic equals $2$. By Poincar\'e Theorem, the sum 
of indices of all singular points of the direction field on $\A$ equals $2$.  
\end{proof}

In fact, for an immersed surface in general position in $\RP^3$ ($\R^3$) 
Theorem~\ref{separates} implies: 

\begin{theorem}
For each connected component of the hyperbolic domain, whose boundary is 
contained in the parabolic curve,  
the flecnodal curve is the union of closed curves each of them 
having an even number (possibly zero) of godrons {\rm (that is, 
of contact points with the boundary of that domain)}. The godrons decompose 
these closed curves into left and right segments.   
\end{theorem}

\begin{corollary}\label{thom}
The boundary of each connected component of the hyperbolic domain of a generic surface 
has an even number of godrons. 
\end{corollary}

The statement of Corollary~\ref{thom} belongs to Thom and Banchoff, \cite{b-Thom}. 
Unfortunately their proof is not exact, since it is based in a wrong statement: 
{\em The Euler characteristic of a 
connected component $H$ of the hyperbolic domain equals the number of godrons in 
$\D H$ at which the hyperbolic domain is locally convex {\rm (the asymptotic line has 
contact with $\D H$ {\em exterior} to $H$)} minus the number of godrons in $\D H$ 
at which the elliptic domain is locally convex {\rm (the asymptotic line has contact 
with $\D H$ {\em interior} to $H$)}}. 
This is wrong: 1) In the bigodron bifurcation (see section~\ref{degenerated_rho=1}) 
two godrons are born (or killed); 2) At both godrons the contact of the 
asymptotic line  with $\D H$ is exterior to $H$ and 3) This bifurcation does not 
change the Euler characteristic of $H$. 

\subsection{Godrons and Swallowtails}\label{fronts}

\noindent 
{\bf Tangential Map and Swallowtails}.\ ~ 
It is well known (c.f. \cite{Salmon}) that under the tangential map of $S$ 
the parabolic curve of $S$ corresponds to the cuspidal edge of $S^\vee$, 
the conodal curve of $S$ corresponds to the self-intersection line of $S^\vee$ 
(this follows from the definitions of dual surface and conodal curve) 
and a godron corresponds to a swallowtail point. 
 

\begin{Lremark}
The most natural approach to the singularities 
of the tangential map is via Arnold's theory of Legendrian 
singularities \cite{avg}. The image of a Legendrian map is called 
the {\em front} of that map. The tangential map of a surface is a Legendre 
map, and so it can be expected to have only Legendre singularities. Thus 
for a surface in general position, the only singularities of its dual surface 
(i.e. of its front) can be: self-intersection lines, cuspidal edges 
and swallowtails. So the godrons are the most complicated singularities of the 
tangential map of a generic surface. 
\end{Lremark}

\noindent 
{\bf Definition of Front}.\ ~ 
In this paper, a {\em front in general position} is a surface $S$ whose 
singularities, and the singularities of its dual surface $S^\vee$, are at most: 
self-intersection lines, semi-cubic 
cuspidal edges and swallowtails. Moreover, we require that the 
parabolic curve never passes through a swallowtail point (the same requirement for the dual 
front). In other words, we are requiring the Legendrian manifold $L_S=L_{S^\vee}$ in 
$PT^*\RP^3=PT^*(\RP^3)^\vee$ (of the contact elements of $\RP^3$ tangent to $S$) to 
be in general position with respect to both natural Legendrian fibrations
$\pi:PT^*\RP^3\ra\RP^3~$ and $~\pi^\vee:PT^*(\RP^3)^\vee\ra(\RP^3)^\vee$. 
Thus, for a front in general position all godrons are simple.

\begin{example*}[the standard swallowtail]
The {\em standard swallowtail} is the discriminant of the vector space of polynomials 
of the form $x^4+ax^2+bx+c$ consisting of all those points $(a,b,c)\in\R^3$ for which 
the polynomial has a multiple root. This discriminant is the tangent developable of the 
curve $\g :t\mapsto (-6t^2, 8t^3, -3t^4)$ consisting of such polynomials having a triple 
real root. So, as any developable surface, it is the envelope of a 1-parameter family of 
planes (the osculating planes of the above curve $\g$, in this case), and hence its 
dual ``surface'' (that is, its image under the tangential map) is just a curve. This 
shows that, from the point of view of projective (or affine) differential geometry, 
the standard swallowtail is absolutely not in general position (in particular, all its points exterior 
to its cuspidal edge $\g$ are parabolic). 
\end{example*}

\noindent 
{\bf The cr-invariant of a swallowtail}.\ ~ 
We can associate a projective invariant (a number) to a swallowtail point $s$
of a generic front $S$: We apply the tangential map of $S$ (in a neighbourhood of $s$) 
to obtain a locally smooth surface $S^\vee$ having a godron with cr-invariant 
$\rho$. The number $\rho(s):=\rho$ is associated to the swallowtail $s$.
\medskip 

\noindent 
{\em The tangential map of $S$ sends the elliptic (hyperbolic) domain of $S$ to the 
elliptic (resp. hyperbolic) domain of $S^\vee$}. Thus the hyperbolic and elliptic domains 
of a front in general position are separated by the cuspidal edge 
(and by the parabolic curve). This implies that there are two types of swallowtails: 

\begin{definition*}
A swallowtail point of a generic front is said to be {\em hyperbolic} ({\em elliptic}) 
if, locally, the self-intersection line of that front is contained in the hyperbolic 
(resp. elliptic) domain. 
\end{definition*} 


The proofs of the following theorems show that the configurations 
of the curves $F$, $P$ and $D$ at a godron have a relevant 
meaning for the local (projective, affine or Euclidean) differential properties 
of the swallowtails.

\begin{theorem}\label{positive-elliptic}
The dual of a surface at a positive godron is an elliptic swallowtail. 
The dual of a surface at a negative godron is a hyperbolic swallowtail. 
\end{theorem}

\begin{proof}
By Proposition~\ref{positive-rho=1}, a godron $g$ is positive (negative) if and only if 
its cr-invariant satisfies $\rho(g) >1$ (resp. $\rho(g)<1$). 

By Theorem~\ref{rho-classification}, 
$\rho(g)>1$ (resp. $\rho(g)<1$) if and only if the conodal curve at $g$ lies locally in 
the elliptic (hyperbolic) domain. 

Finally, since the tangential map sends the elliptic 
(hyperbolic) domain to the elliptic (resp. hyperbolic) domain of the dual surface, 
it is evident that the conodal curve at $g$ lies locally in the elliptic (hyperbolic) 
domain if and only if the dual surface is an elliptic (resp. hyperbolic) swallowtail.  
\end{proof}

\begin{figure}[ht]
\centerline{\psfig{figure=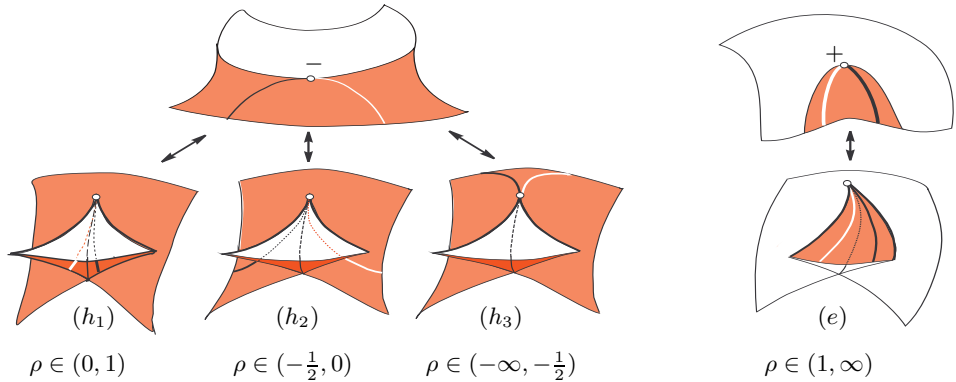,height=4.6cm}}
\begin{picture}(0,0)
\put(29,20){\small $(h_1)$} 
\put(106,20){\small $(h_2)$} 
\put(183,20){\small $(h_3)$} 
\put(311,20){\small $(e)$} 
\put(117,116){\small $\bf{-}$} 
\put(314,120){\small $\bf{+}$}
\put(13,2){\small $\rho\in(0,1)$} 
\put(90,2){\small $\rho\in(-\frac{1}{2},0)$} 
\put(163,2){\small $\rho\in(-\infty,-\frac{1}{2})$} 
\put(291,2){\small $\rho\in(1,\infty)$} 
\end{picture}
\caption{\small The $4$ generic configurations of the flecnodal curve at a swallowtail.}
\label{duality}
\end{figure}

\begin{theorem}\label{4-fronts}
In the neighbourhood of a swallowtail point $s$ of a front $S$ in general position, 
the flecnodal curve $F$ has a cusp whose tangent direction coincides with that of 
the cuspidal edge. The point $s$ separates $F$ locally into its left and 
right branches. There are four possible generic configurations of $F$ in the 
neighbourhood of $s$ {\rm (see Fig.~\ref{duality}):}
\smallskip

\noindent
$(e)$ For an elliptic swallowtail the flecnodal curve is a 
cusp lying in the small domain bounded by the cuspidal edge 
{\rm ($\rho(s)\in(1,\infty)$)}. 
\smallskip

\noindent
There are $3$ different generic types of hyperbolic swallowtails.
\smallskip  

\noindent
$(h_1)$ Each branch of the cuspidal edge is separated from the 
self-intersection line by one branch of the flecnodal curve 
{\rm ($\rho(s)\in(0,1)$)}. 

\noindent
$(h_2)$ The self-intersection line lies between the two branches 
of the flecnodal curve and separates them from the branches of the 
cuspidal edge. The cusp of the flecnodal curve 
points in the same direction as the cusp of the cuspidal edge 
{\rm ($\rho(s)\in(-\frac{1}{2},0)$)}.

\noindent
$(h_3)$ The cusp of the flecnodal curve and the cusp of the cuspidal edge 
are pointing in opposite directions {\rm ($\rho(s)\in(-\infty,-\frac{1}{2})$)}.
\end{theorem}

\subsection{Finer classification of swallowtails}\label{finer_swallowtail}
Besides the description of the $4$ generic configurations of the flecnodal curve 
in the neighbourhood of a swallowtail point (given in Theorem~\ref{4-fronts} and 
Fig.~\ref{duality}), it is interesting to know how the cuspidal edge, 
the flecnodal curve and the self-intersection line of the swallowtail surface are 
placed with respect to the tangent plane. 

The following theorem is a refinement of Theorem~\ref{4-fronts}, providing the 
local configurations of the cuspidal edge, the flecnodal curve and the self-intersection 
line with respect to the tangent plane, at a swallowtail.

\begin{figure}[ht]
\centerline{\psfig{figure=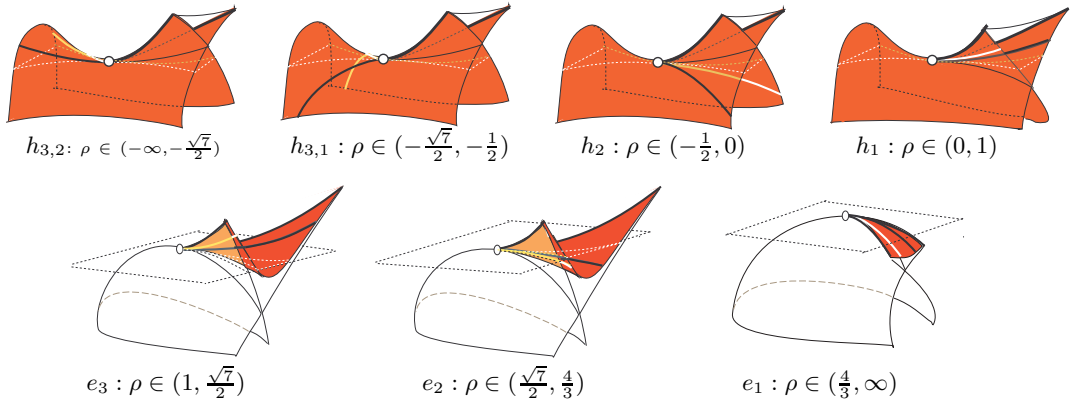,height=4.8cm}}
\begin{picture}(0,0)
\put(-10,92){\small $h_{3,2}\scriptstyle{:\ \rho\ \in\ (-\infty,-\frac{\sqrt{7}}{2})}$} 
\put(90,92){\small $h_{3,1}:\rho\in(-\frac{\sqrt{7}}{2},-\frac{1}{2})$} 
\put(200,92){\small $h_2: \rho\in(-\frac{1}{2},0)$} 
\put(304,92){\small $h_1:\rho\in(0,1)$} 
\put(13,2){\small $e_3:\rho\in(1,\frac{\sqrt{7}}{2})$} 
\put(140,2){\small $e_2:\rho\in(\frac{\sqrt{7}}{2},\frac{4}{3})$} 
\put(261,2){\small $e_1:\rho\in(\frac{4}{3},\infty)$} 
\end{picture}
\caption{\small The 7 generic configurations of the flecnodal curve, the self-intersection 
line, the cuspidal edge and the tangent plane at a swallowtail.}
\label{7_fronts_fig}
\end{figure}

\begin{theorem}\label{7-fronts}
In the notations -- and conclusions -- of Theorem~\ref{4-fronts}, write $\Sigma$
for the cuspidal edge of the front $S$ and $D^\vee$ for its self-intersection line. \\
If $s$ is elliptic and $1<\rho<4/3$ (resp. hyperbolic, $\rho<1$), then the tangent plane to $S$ at 
$s$ intersects $S$ along two semi-cubic cusps $C_-^\vee$ and $C_+^\vee$, pointing 
in the same direction ---as the cuspidal edge-- (resp. in opposite directions) and 
lying in the elliptic and in the hyperbolic domain, respectively (resp. lying both
in the hyperbolic domain). There are $7$ possible generic configurations of $F$, $\Sigma$, 
and $D^\vee$ with respect to the tangent plane {\rm (see Fig.~\ref{7_fronts_fig})}. 
\medskip

\noindent
There are $3$ different generic configurations for elliptic swallowtails.
\medskip

\noindent
$(e_1)$ The surface lies locally on one side of the tangent plane at $s$ 
$(\rho>4/3)$. \\
$(e_2)$ The cusp $C_+^\vee$ separates locally the cuspidal edge $\Sigma$ from  
the flecnodal curve $F$, and $C_-^\vee$ separates $\Sigma$ from the self-intersection 
line $D^\vee$ $(\rho\in(\frac{\sqrt{7}}{2},\frac{4}{3}))$. \\
$(e_3)$ The flecnodal curve $F$ separates locally $\Sigma$ from $C_+^\vee$, and again
$C_-^\vee$ separates $\Sigma$ from $D^\vee$ $(\rho\in(1,\frac{\sqrt{7}}{2}))$.
\medskip

\noindent
There are $4$ different generic configurations for hyperbolic swallowtails.
\medskip

\noindent
$(h_1)$ The cusp $C_+^\vee$ separates locally $F$ from $D^\vee$ 
$(\rho\in(0,1))$. \\
$(h_2)$ Again, the cusp $C_+^\vee$ separates locally $F$ from $D^\vee$ 
$(\rho\in(-\frac{1}{2},0))$. \\
$(h_{3,1})$ The tangent plane separates locally $F$ from $\Sigma$ and $D^\vee$
$(\rho\in(-\frac{\sqrt{7}}{2},-\frac{1}{2}))$.\\
$(h_{3,2})$ The curves $F$, $\Sigma$ and $D^\vee$ lie locally on the same side of the 
tangent plane to $S$ at $s$ $(\rho<-\frac{\sqrt{7}}{2})$. 
\end{theorem}

\subsubsection{The local transition at swallowtails with $\rho=0$}\label{flec-swallowtail-section}

Now we can describe the local bifurcation of the tangential singularities 
occurring in a generic one parameter family of wave fronts at the 
moment of a swallowtail with cr-invariant $\rho=0$. 

Let $S_t$, $t\in\R$, a generic one parameter family of fronts 
such that for $t=0$ the front $S_0$ has a swallowtail point $s_0$ 
with $\rho=0$ (we name it a {\em flec-swallowtail}). 

Just by duality of the flec-godron transition, we have that as $t$ is 
increasing and passing through $0$, a biflecnode $b_t$ of $S_t$ 
is moving along the flecnodal curve, passing, at $s_0$, from one branch to the other of 
the flecnodal curve. So, `a flec-swallowtail is the superposition of a biflecnode and 
a swallowtail'.
 
\begin{figure}[ht]
\begin{picture}(0,0)
\put(58,0){\small $\rho<0$}
\put(159,0){\small $\rho=0$} 
\put(257,0){\small $\rho>0$}
\end{picture}
\centerline{\psfig{figure=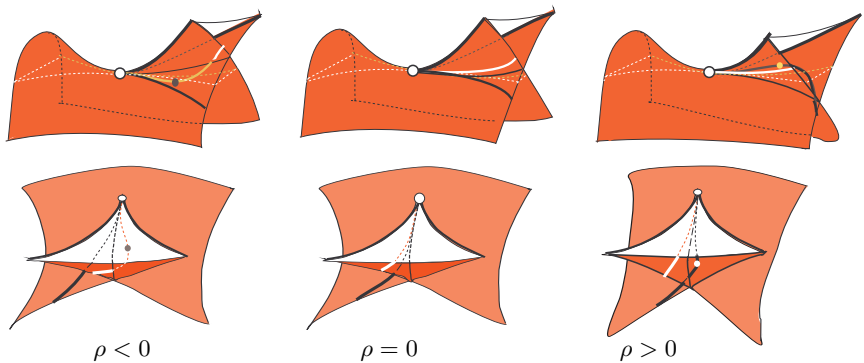,height=4.7cm}}
\caption{\small Two points of view of the transition at a swallowtail with $\rho=0$ 
}
\label{flex-swallowtail}
\end{figure}

In Fig.~\ref{flex-swallowtail}, we show the pictures describing the bifurcation 
of the flecnodal curve occurring in a generic one parameter family of wave fronts at 
a flec-swallowtail. They follow from \S \ref{degenerated_rho=0}, 
Theorem~\ref{7-fronts}, and Figures~\ref{flexgodron} and \ref{7_fronts_fig}.

At the moment with $\rho=0$ both the self-intersection line $D^\vee$ and 
the flecnodal curve $F$ have higher order of contact with the tangent plane; 
in particular, for $\rho>0$, $D^\vee$ lies locally in the opposite side 
of the tangent plane than for $\rho<0$. The flecnodal curve is also locally passing 
from one side to the other of the tangent plane at $\rho=0$. 


\subsection{The local $q$-contour of a surface}\label{g_contour}

We shall investigate  whether any arbitrarily small neighbourhood of a godron $g$ 
(or of a swallowtail point $s$) of a front $S$ in general position has (or not) points $p\in S$ 
such that the tangent plane to $S$ at $p$ passes through $g$ (resp. through $s$). 
\smallskip

Consider a point $q$ of a smooth surface or of a wave front $S$ 
in the $3$-space (Euclidean, affine or projective).  

\begin{definition*}
The {\em $q$-contour of $S$} is the set of points of 
tangency of the planes tangent to $S$ passing through $q$. 
\end{definition*}

Equivalently, the {\em $q$-contour of $S$} is the set of points of tangency of the  
lines tangent to $S$ passing through $q$. \ 
The following lemma is evident: 

\begin{lemma*}
The $q$-contour of $S$, without the point $q$, consist of the critical points of the 
``stereographic'' projection of $S$, from $q$ 
to $\RP^2_q$,  $\pi_q:S\smallsetminus \{q\}\ra \RP^2_q$,  {\rm (associating to each 
point $p\in S$ the line joining it to $q$)}.
\end{lemma*}

One usually considers the projection of a surface from a point exterior to it (for which 
the set of critical points is, generically, a -- possibly empty -- smooth curve), but here 
we shall consider the projection of a surface from a point belonging 
to it, and the critical points of such a projection. 

Hence, given a point $q$ of a smooth surface (or of a front) $S$ in the $3$-space, 
it is interesting to know whether any arbitrarily small neighbourhood of $q$ 
contain points of the $q$-contour of $S$ different from $q$ and, if it is the case, 
to know its behaviour at $q$. 

\begin{remark*}
For any point $q$ of a smooth surface (or of a wave front) $S$ the 
$q$-contour of $S$ is projectively invariant. 
\end{remark*}

\begin{definition*}
The {\em local $q$-contour} of a smooth surface (or of a front) $S$ 
is the germ at $q$ of the $q$-contour of $S$. We denote it by $C^q(S)$. 
\end{definition*}

We will say that the local $q$-contour of a smooth surface (or of a front) is trivial if it
consists just of the point $q$. 

\begin{example*}
For the elliptic points of a surface (or of a front) $S$ the $q$-contour of $S$ is trivial. \ 
For the hyperbolic points the $q$-contour of $S$ consist of two transverse curves (each one 
being tangent to one of the asymptotic lines). 
\end{example*}

The following theorem shows, for example, that there are two essentially different kinds of positive 
godrons: For one kind the local $g$-contour is non-trivial and for the other it is trivial. 
This difference is more visible on the dual surface: The tangent plane traverses locally
the surface or not, (see Fig.\ref{7_fronts_fig}). 

\begin{theorem}\label{4/3}
Let $g$ be a godron of a smooth surface $S$.

\noindent
$(a)$ There exists a neighbourhood $U$ of $g$, in $S$, which contains no point of the set 
$C^g(S)\smallsetminus \{g\}$ {\rm (i.e. no tangent plane to 
$U\smallsetminus \{g\}$ passes through $g$)} if and only if $\rho(g)>4/3$ ; 
\medskip

\noindent
$(b)$ If $\rho(g)<4/3$, $\rho(g)\neq 1$, then, in any sufficiently small neighbourhood of $g$, the 
$g$-contour of $S$ consists of two smooth curves tangent to the parabolic curve at $g$. 
If the godron $g$ is positive $(1<\rho<4/3)$, then one of these two curves lies locally in the 
elliptic domain and the other in the hyperbolic domain. If $g$ is negative $(\rho<1)$, 
then both curves lie in the hyperbolic domain. 
\end{theorem}

In case $b$ of Theorem~\ref{4/3}, we write $C_-$ and $C_+$ for the tangent curves 
forming the local $g$-contour of $S$. 
Since $C_-$ and $C_+$ are tangent to the parabolic curve at $g$, 
their projections to the tangent plane can be written locally as $y=c_\sC^-x^2+\ldots$ and 
$y=c_\sC^+x^2+\ldots$ (using Platonova's normal form). The relative positions of $F$, $P$ and $D$ 
with respect to $C_-$ and $C_+$ are determined by the order, in the real line, of the 
coefficients $c_\sF$, $c_\sP$, $c_\sD$, $c_\s=1$, $c_\sC^-$ and $c_\sC^+$.

\begin{figure}[ht]
\centerline{\psfig{figure=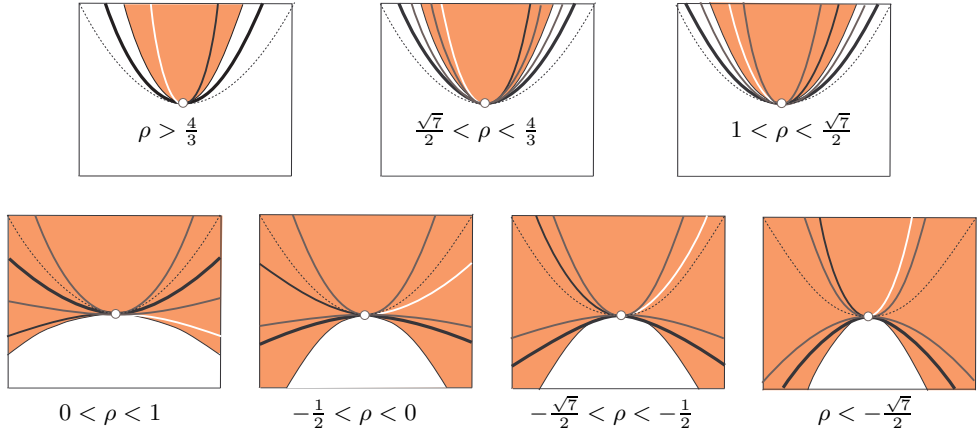,height=5.2cm}}

\begin{picture}(0,0)
\put(50,110){\small $\rho>\frac{4}{3}$} 
\put(154,110){\small $\frac{\sqrt{7}}{2}<\rho <\frac{4}{3}$} 
\put(274,110){\small $1<\rho<\frac{\sqrt{7}}{2}$} 
\put(20,3){\small $0<\rho<1$} 
\put(108,3){\small $-\frac{1}{2}<\rho<0$} 
\put(198,3){\small $-\frac{\sqrt{7}}{2}<\rho<-\frac{1}{2}$} 
\put(307,3){\small $\rho<-\frac{\sqrt{7}}{2}$} 
\end{picture}
\caption{\small The 7 generic configurations of the curves $F$ (half-white half-black curves), 
$P$ (boundary between white and gray domains), $D$ (thick curves), the separating 
$2$-jet (doted curves) and the $g$-contour (gray curves) at a godron $g$.}
\label{classification_contour}
\end{figure}

\begin{theorem}\label{rho-classification_contour}
Given a godron $g$ of a generic smooth surface,   
there are $7$ possible configurations of the curves $F$, $P$ and $D$ 
with respect to the $g$-contour the surface at $g$ 
{\rm (they are represented in Fig.~\ref{classification_contour})}. The actual configuration 
at $g$ depends on which of the following $7$ open intervals 
the cr-invariant $\rho(g)$ belongs to, respectively: 
$$
\begin{array}{lcc}

\rho\in(\frac{4}{3},\infty) &  \iff &  1<c_\sD<c_\sP<c_\sF;\\ 
\rho\in(\frac{\sqrt{7}}{2},\frac{4}{3}) &  \iff &  1<c_\sD<c_\sC^-<c_\sP<c_\sC^+<c_\sF;\\ 
\rho\in(1,\frac{\sqrt{7}}{2}) &  \iff &  1<c_\sD<c_\sC^-<c_\sP<c_\sF<c_\sC^+;\\
\rho\in(0,1) &  \iff &  c_\sP<c_\sF<c_\sC^-<c_\sD<1<c_\sC^+;\\
\rho\in(-\frac{1}{2},0) &  \iff &  c_\sP<c_\sD<c_\sC^-<c_\sF<1<c_\sC^+;\\
\rho\in(-\frac{\sqrt{7}}{2},-\frac{1}{2}) &  \iff &  c_\sP<c_\sD<c_\sC^-<1<c_\sF<c_\sC^+;\\
\rho\in(-\infty,-\frac{\sqrt{7}}{2}) &  \iff &  c_\sP<c_\sD<c_\sC^-<1<c_\sC^+<c_\sF.
\end{array}
$$
\end{theorem}

\subsubsection{The local $q$-contour at a swallowtail point}\label{g_contour}

It is also interesting to know whether the local $s$-contour at a swallowtail point $s$
of a front in general position is trivial or not, and, if it is non trivial, to know how it is
placed with respect to the self-intersection line $D^\vee$, the cuspidal edge $\Sigma$ and 
the flecnodal curve $F$. \ 
The following theorem is a refinement of Theorem~\ref{4-fronts}, providing 
the required local configurations.

\begin{figure}[ht]
\centerline{\psfig{figure=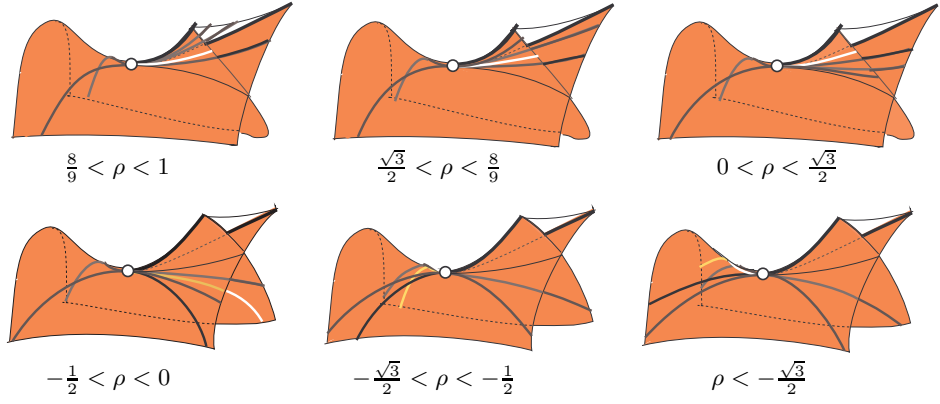,height=4.7cm}}
\begin{picture}(0,0)
\put(29,82){\small $\frac{8}{9}<\rho<1$} 
\put(147,82){\small $\frac{\sqrt{3}}{2}<\rho <\frac{8}{9}$} 
\put(276,82){\small $0<\rho<\frac{\sqrt{3}}{2}$} 
\put(22,2){\small $-\frac{1}{2}<\rho<0$} 
\put(138,2){\small $-\frac{\sqrt{3}}{2}<\rho<-\frac{1}{2}$} 
\put(274,2){\small $\rho<-\frac{\sqrt{3}}{2}$} 
\end{picture}
\caption{\small The 7 generic configurations of the flecnodal curve, the self-intersection 
line, the cuspidal edge and the local $s$-contour (gray cusps) at a swallowtail $s$.}
\label{7_fronts_contour}
\end{figure}

\begin{theorem}\label{6-fronts}
In the notations -- and conclusions -- of Theorem~\ref{4-fronts}\:  

\noindent
$(e)$ The local $s$-contour of $S$ is trivial {\rm (consisting just of $s$)} 
if and only if $s$ is an elliptic swallowtail ($\rho>1$). \\
$(h)$ If $s$ is hyperbolic, then the local $s$-contour of $S$ consists of two 
cusps $T_-^\vee$, $T_+^\vee$, pointing in opposite directions. 
In this case, there are $6$ possible generic configurations 
of the local $s$-contour with respect to $F$, $\Sigma$ and $D^\vee$ 
{\rm (see Fig.~\ref{7_fronts_contour})}{\rm :} \\
$(h_{1,1})$ $T_-^\vee$ lies in the elliptic domain $(\rho\in(\frac{8}{9},1))$. \\
$(h_{1,2})$ $T_-^\vee$ separates locally $\Sigma$ from $F$
            $(\rho\in(\frac{\sqrt{3}}{2}, \frac{8}{9}))$. \\
$(h_{1,3})$ $T_-^\vee$ separates locally $F$ from $D^\vee$
            $(\rho\in(0, \frac{\sqrt{3}}{2}))$. \\
$(h_2)$ $T_-^\vee$ separates locally $D^\vee$ from $F$ $(\rho\in(-\frac{1}{2},0))$. \\
$(h_{3,1})$ $T_-^\vee$ separates $D^\vee$ from $F$. $T_-^\vee$ and $T_+^\vee$ 
are separated by $F$ $(\rho\in(-\frac{\sqrt{3}}{2},-\frac{1}{2}))$. \\
$(h_{3,2})$ $T_-^\vee$ and $F$ are separated locally by $T_+^\vee$ 
$(\rho\in(-\infty,-\frac{\sqrt{3}}{2}))$.
\end{theorem}

\subsection{The tangent section of a surface or of a front}\label{t_section}

The intersection of a surface $S$ with its tangent plane at a point $q$ 
will be called the {\em $q$-tangent section} of $S$ or the {\em tangent section}
of $S$ at $q$. 

\begin{definition*}
The {\em local tangent section} of a smooth surface (or of a front) $S$ 
at its point $q$ 
is the germ at $q$ of the $q$-tangent section of $S$. The local tangent 
section at $q$ is said to be {\em trivial} if it consists just of the point $q$. 
\end{definition*}

\begin{lemma}\label{tangent_section}
Consider a simple godron $g$ of a smooth surface $S$. 
\medskip

\noindent
$(a)$ The local tangent section of $S$ at $g$ is trivial (consisting just of $g$) 
if and only if $g$ is a positive godron {\rm ($\rho>1$)}. 
\medskip

\noindent
$(b)$ If $g$ is negative {\rm ($\rho<1$)}, then the local tangent section of $S$ at $g$ 
consists of two smooth curves (simply) tangent to the parabolic curve at $g$. 
\end{lemma}

Write $T_-$ and $T_+$ for the tangent curves of the case $b$ of Lemma~\ref{tangent_section},  
forming the local tangent section of $S$ at $g$. 
Since $T_-$ and $T_+$ are tangent to the parabolic curve at $g$, 
their projections to the tangent plane can be written locally as $y=c_\sT^-x^2+\ldots$ and 
$y=c_\sT^+x^2+\ldots$ (using Platonova's normal form). The relative positions of $F$, $P$ and $D$ 
with respect to $T_-$ and $T_+$ are determined by the order, in the real line, of the 
coefficients $c_\sF$, $c_\sP$, $c_\sD$, $c_\s=1$, $c_\sT^-$ and $c_\sT^+$.

\begin{figure}[ht]
\centerline{\psfig{figure=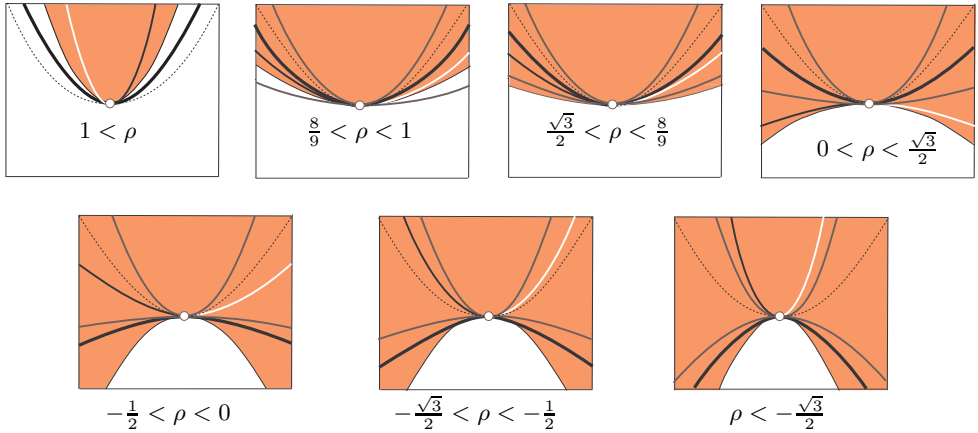,height=5.2cm}}
\begin{picture}(0,0)
\put(28,110){\small $1<\rho$} 
\put(114,110){\small $\frac{8}{9}<\rho<1$} 
\put(204,110){\small $\frac{\sqrt{3}}{2}<\rho<\frac{8}{9}$} 
\put(307,103){\small $0<\rho<\frac{\sqrt{3}}{2}$} 
\put(38,3){\small $-\frac{1}{2}<\rho<0$} 
\put(147,3){\small $-\frac{\sqrt{3}}{2}<\rho <-\frac{1}{2}$} 
\put(274,3){\small $\rho<-\frac{\sqrt{3}}{2}$} 
\end{picture}
\caption{\small The 7 generic configurations of the curves $F$ (half-white half-black curves), 
$P$ (boundary between white and gray domains), $D$ (thick curves), the separating 
$2$-jet (doted curves) and the local tangent section (gray curves) at a godron.}
\label{classification_section}
\end{figure}

\begin{theorem}\label{rho-classification_tangent}
Given a simple godron $g$ of a smooth surface,   
there are $7$ possible configurations of the curves $F$, $P$ and $D$ 
with respect to the separating $2$-jet and the local tangent section at $g$ 
{\rm (they are represented in Fig.~\ref{classification_section})}. The actual configuration 
at $g$ depends on which of the following $7$ open intervals 
the cr-invariant $\rho(g)$ belongs to, respectively\: 
$$
\begin{array}{lcl}
\rho\in(1,\infty) &  \iff &  \ ~ 1~<c_\sD<c_\sP<c_\sF\ \mbox{ \rm ({\small $T^g(S)$ is trivial})};\\ 
\rho\in(\frac{8}{9},1) &  \iff &  c_\sT^-<c_\sP<c_\sF<c_\sD<1<c_\sT^+;\\
\rho\in(\frac{\sqrt{3}}{2},\frac{8}{9}) &  \iff &  c_\sP<c_\sT^-<c_\sF<c_\sD<1<c_\sT^+;\\
\rho\in(0,\frac{\sqrt{3}}{2})  & \iff &  c_\sP<c_\sF<c_\sT^-<c_\sD<1<c_\sT^+;\\
\rho\in(-\frac{1}{2},0) &  \iff &  c_\sP<c_\sD<c_\sT^-<c_\sF<1<c_\sT^+;\\
\rho\in(-\frac{\sqrt{3}}{2},-\frac{1}{2}) &  \iff &  c_\sP<c_\sD<c_\sT^-<1<c_\sF<c_\sT^+;\\
\rho\in(-\infty,-\frac{\sqrt{3}}{2}) &  \iff &  c_\sP<c_\sD<c_\sT^-<1<c_\sT^+<c_\sF.
\end{array}
$$
\end{theorem}

\noindent
{\small {\bf Note}. By Theorem~\ref{rho-classification} (Fig.~\ref{classification}) the curves $P$ and $F$ 
change their local convexity at the values $\rho=\frac{2}{3},\ \rho=\frac{1}{2}$
in the interval $0<\rho<\frac{\sqrt{3}}{2}$, respectively. So, in this interval 
we can choose $3$ pictures with the same configuration but with different local convexities of 
$P$ and $F$. In Fig.~\ref{classification_section} we have chosen the picture corresponding to 
the subinterval $(0,\frac{1}{2})\subset(0,\frac{\sqrt{3}}{2})$. In Fig.~\ref{classification_contour} 
we had the same situation and we chosen the same subinterval $(0,\frac{1}{2})\subset(0,1)$.}

\begin{remark*}
The above theorems show that the value of the cr-invariant of a 
simple godron $g$ of a smooth surface $S$ determines completely the 
local configuration of the parabolic, flecnodal and conodal curves with respect to 
the local tangent section at $g$ and to the local $g$-contour of $S$ (determining 
also their local convexities). Similarly, the local configurations obtained for the swallowtails 
are also determined by the cr-invariant. 
\end{remark*}

\begin{remark*}
In \cite{B_G_T_95}, it was observed that there are two types of godrons, 
called {\em elliptic} and {\em hyperbolic}, corresponding to the sign of 
a `discriminant' relating some coefficients of the $4$-jet of the surface. 
In fact, they are the positive and negative godrons of \S \ref{index-definition} 
(definition of index), which are distinguished geometrically by the behaviour of 
the asymptotic directions along the parabolic curve (see Fig.~\ref{index}). 
Here, we have provided several geometric characterisations of the 
godrons of indices  $+1$ and $-1$: In terms of the cr-invariant
(Proposition~\ref{positive-rho=1}), in terms of the left and right flecnodal curves 
(Theorem~\ref{lr}), in terms of the geometry of the dual surface -- elliptic and 
hyperbolic swallowtails (Theorem~\ref{positive-elliptic}) 
and in terms of the asymptotic double (Figure~\ref{double}). Of course, 
Lemma~\ref{tangent_section} also distinguishes these points (in terms of the 
local tangent section). 
\end{remark*}

\section{The proofs of the theorems}\label{proofs}
{\bf Preparatory conventions and results}.\ ~ 
In the sequel, we will consider the surface $S$ as the graph 
of a smooth function $z=f(x,y)$, where $x,y,z$ form an affine coordinate 
system. The asymptotic directions satisfy the equation:
$$f_{xx}(dx)^2+2f_{xy}dxdy+f_{yy}(dy)^2=0.$$ 
For $dy=pdx$, this equation takes the form 
$$A^f(x,y,p)=f_{xx}+2f_{xy}p+f_{yy}p^2=0.\eqno(1)$$
Equation (1) is called the {\em asymptote-equation} of $f$. 
\medskip

In what follows, we will assume without loss of generality that the point 
under consideration in the $(x,y,p)$-space is the origin: by a translation 
and a rotation in the $(x,y)$-plane, 
we can take $(x,y)=(0,0)$ and $p=0$, respectively. 

Moreover, we will take an affine coordinate system $x,y,z$ such that the 
$(x,y)$-plane is tangent to $S$ at the point under 
consideration. Thus we will have the conditions 
$$f(0,0)=f_x(0,0)=f_y(0,0)=0.\eqno(2)$$ 

The parabolic curve of the surface $z=f(x,y)$ is the restriction 
of the graph of $f$ to the discriminant curve (in the $(x,y)$-plane) 
of equation (1). That is, the parabolic curve is determined by the equations 
$$A^f(x,y,p)=0 \ \mbox{ and }\ A^f_p(x,y,p)=0.\eqno(*)$$

The fact that a godron is a folded singularity of eq.~(1) implies that 
$$A^f_x(0,0,0)=0.\eqno(**)$$

The conditions $(*)$ and $(**)$, at the origin in the $(x,y,p)$-space, imply that 
$$f_{xx}=f_{xy}=f_{xxx}=0\eqno(3)$$ 
at the origin in the $(x,y)$-plane.
\medskip 

The choice of a coordinate system such that the $x$-axis is an asymptotic 
direction of $S$ at the origin is equivalent to our assumption that the point 
under consideration in the $(x,y,p)$-space is the origin. 

So the $x$-axis is tangent to the parabolic curve at the godron.

\subsection{Proof of Theorem~\ref{godron_flattening}}\label{proof2}

Let $\g(t)=(x(t),y(t),z(t))$ be a curve on $S$, where $z(t)=f(x(t),y(t))$,  
which is tangent to the parabolic curve at the origin, that is, 
$$\dot{y}(0)=0.\eqno(4)$$ 

Since all our calculations and considerations take place at 
the origin $(x,y)=(0,0)$ and at $t=0$, we will omit to write this explicitly. 

Evidently conditions $(2)$ imply $\dot{z}=f_x \dot{x}+f_y \dot{y}=0$. 
The equality
$$\ddot{z}=f_x \ddot{x}+f_y \ddot{y}+(f_{xx}\dot{x}^2+2f_{xy}\dot{x}\dot{y}+f_{yy}\dot{y}^2)$$
together with conditions (2), (3) and (4) imply that $\ddot{z}=0$. This proves that 
the plane $z=0$ is osculating. 

Finally, the equality 
$$
\begin{array}{rcr}
\dddot{z} & = & f_x \dddot{x}+f_y \dddot{y}+3(f_{xx}\dot{x}\ddot{x}+
f_{xy}(\dot{x}\ddot{y}+\ddot{x}\dot{y})+f_{yy}\dot{y}\ddot{y})  \\
 &  & \\
  &  & +\ f_{xxx}\dot{x}^3+3f_{xxy}\dot{x}^2\dot{y}+3f_{xyy}\dot{x}\dot{y}^2+f_{yyy}\dot{y}^3
\end{array}$$
together with conditions (2), (3) and (4), imply that $\dddot{z}=0$, proving that 
the first three derivatives of $\g$ at $t=0$ are linearly dependent (all of them 
lie in the $(x,y)$-plane). So $\g$ has a flattening or an inflection at the origin, 
according to the linear independence or dependence, respectively, 
of its first two derivatives at $t=0$. \fin

\def\F{\bar{F}}
\def\P{\bar{P}}
\def\D{\bar{D}}
\subsection{Preliminary remarks and computations}\label{remarks-computations}

We recall that Platonova's Theorem \cite{Platonova} implies 
that at a godron of a generic smooth surface $S$, there is an 
affine coordinate system such that $S$ is locally given by 
$$z=\frac{y^2}{2}-x^2y+\lambda x^4+\f(x,y) \ \ \ 
\mbox{(for some $\lambda\neq \frac{1}{2},0$)} \eqno(G1)$$
where $\f$ is the sum of homogeneous polynomials in $x$ and $y$ of degree 
greater than $4$ and (possibly) of flat functions. 

The information we need about $S$ (for the proofs of our theorems) is contained 
in its $4$-jet. The fourth degree terms are enough to find the generic configurations 
of the curves $F$, $P$, and $D$ (with respect to the other special curves) and to 
identify the exceptional values of $\rho$, separating the different generic 
configurations (it is not difficult to see, for instance, that the 
canonical coefficients of the curves $F$, $P$ and $D$ are independent of the term $\f$ in (G1)).
Thus, in the proofs of our theorems, 
we will systematically use Platonova's normal form of the $4$-jet of $S$. 
The terms of degree $5$ can be relevant, however, to study the bifurcations 
occurring at the exceptional values of $\rho$ (for example, $\rho=0$), separating 
the values corresponding to the generic configurations. 
Such bifurcations will be studied in another paper. 
\medskip

First we need to calculate the curves $F$, $P$ and $D$. For we 
need the second partial derivatives of the functions 
$f(x,y;\lambda)=\frac{y^2}{2}-x^2y+\lambda x^4$:  
$$f_{xx}=-2y+12\lambda x^2,\ \ \ f_{xy}=-2x,\ \ \ f_{yy}=1.\eqno(H)$$ 

The asymptote-equations of the surfaces $z=\frac{y^2}{2}-x^2y+\lambda x^4$ 
are therefore given by 
$$A^f(x,y,p;\lambda)=(12\lambda x^2-2y)-4xp+p^2=0.\eqno(5)$$ 

We are interested in the configurations of the curves 
$F$, $P$ and $D$, at the godron $g$. According to Theorem 2, these curves 
have at least $4$-point contact with the $(x,y)$-plane. We will thus consider 
the curves $\F$, $\P$ and $\D$, on the $(x,y)$-plane, whose images by $f$ are 
$F$, $P$ and $D$, respectively. These plane curves have the same $2$-jet as 
$F$, $P$ and $D$, respectively. 
\medskip 

\noindent 
{\bf The parabolic curve}.\ ~ 
The equations $(*)$ of \S\ref{proof2} imply that $\P$ is given by the 
Hessian of $f$, $f_{xy}^2-f_{xx}f_{yy}=0$. From (H), one obtains 
that $\P$ is a parabola: 
$$y=2(3\lambda -1)x^2.$$ 

\noindent
{\bf The flecnodal curve}.\ ~
According to \cite{Uribetesis,surf-evolution}, the curve $\F$ 
associated to the surface $z=f(x,y)$ is obtained from the intersection of the surfaces 
$$A^f(x,y,p)=0\ \ \mbox{and }\ \ I^{A^f}(x,y,p):=(A^f_x+pA^f_y)(x,y,p)=0,$$ 
in the $(x,y,p)$-space, by the projection of this intersection to the 
$(x,y)$-plane, along the $p$-direction. From eq.~(5) one obtains 
$$I^{A^f}(x,y,p)=6(4\lambda x-p).$$
Combining the equation $p=4\lambda x$ with eq.~(5) one obtains that $\F$ 
is a parabola: 
$$y=2\lambda (4\lambda-1)x^2.$$

\noindent
{\bf The conodal curve}.\ ~
Since Platonova's normal form is symmetric with respect to the $x$-direction, 
the bitangent planes in the neighbourhood of $g$ are invariant under the reflection 
$(x,y,z)\mapsto (-x,y,z)$. Thus the points of the conodal curve satisfy 
$f_x(x,y;\lambda)=0$. That is, $-2x(y-2\lambda x^2)=0$. Thus the curve $\D$ is a 
parabola: 
$$y=2\lambda x^2.$$

\subsection{Proof of Theorem \ref{rho=2l}}
We consider the parabolas $\F$, $\P$ and $\D$ as graphs of functions $y=y(x)$.
The Legendrian curves $L_F$, $L_P$ and $L_D$ in the $(x,y,p)$-space $J^1(\R,\R)$ 
(which is the space of $1$-jets of the real functions $y(x)$ of 
one real variable) are tangent to the contact plane 
$\Pi$ at the origin (parallel to the plane $y=0$). 
The slope of the tangent line at the origin, of each of these Legendrian curves, 
equals twice the second derivative at zero of the function $y=y(x)$ 
associated to the corresponding parabola, that is, equals twice the coefficient 
of that parabola (note that the term $\f$ in $(G1)$ will contribute with higher order 
terms which will have no influence on these coefficients). 

The Legendrian curve consisting of the contact elements tangent to the origin 
is vertical. Write $\ell_g$ for its tangent line. The cross-ratio of the tangent 
lines $\ell_F$, $\ell_P$, $\ell_D$ and $\ell_g$ is given in terms of the coefficients $c$ 
of the parabolas $\F$, $\P$ and $\D$ by 

$$\rho(g)=(\ell_F,\ell_P,\ell_D,\ell_g)= \frac{c(F)-c(D)}{c(P)-c(D)}
=\frac{2\lambda (4\lambda-1)-2\lambda}{2(3\lambda-1)-2\lambda}=2\lambda.$$ 
This proves Theorem~\ref{rho=2l}. \fin
\bigskip


\noindent
{\bf Rewriting the equations in terms of $\rho$}.\ ~
After Theorem \ref{rho=2l}, we rewrite Platonova's normal forms of the $4$-jet 
of $S$ at a godron 
and the equations of the curves $\F$, $\P$ and $\D$ in terms of the cr-invariant $\rho$: 
$$z=\frac{y^2}{2}-x^2y+\rho\frac{x^4}{2} \ \ \ \ (\rho \neq 1,0). \eqno(R)$$
$$y=(3\rho -2)x^2; \eqno(P)$$
$$y=\rho(2\rho-1)x^2;\eqno(F)$$
$$y=\rho x^2.\eqno(D)$$

\subsection{Proof of the Separating $2$-jet Lemma}\label{proof_separating}

An easy way to compute (and to see) 
the dual surface of $S\subset \R^3$, viewed as a surface in 
the same space $\R^3$ and with the same coordinate system, 
is by the `polar duality map' with respect to a quadric. 
The calculations are simpler if the quadric (considered
for this map) is a paraboloid of revolution (see \cite{Uribepdll}). 
Moreover, if the surface $S$ is the graph of a function $z=f(x,y)$, 
then the polar duality map with respect to the paraboloid $z=\frac{1}{2}(x^2+y^2)$ 
coincides with the classical Legendre transform of $f$. 
So, the dual surface of the graph $\{ (x,y,f(x,y))\}$ has the following parametrisation: 
$$\tau_f:(x,y)\mapsto (f_x(x,y),\ f_y(x,y),\ xf_x(x,y)+yf_y(x,y)-f(x,y)).$$
In the case of the surfaces $S_\rho$ given in eq.~$(R)$, one obtains 
$$\tau_\rho:(x,y)\mapsto 
\left(-2xy+2\rho x^3,\  y-x^2,\  \frac{y^2}{2}-2x^2y+3\rho \frac{x^4}{2}\right).\eqno(R^\vee)$$

The images of our plane curves $\F$, $\P$ and $\D$, under $\tau_\rho$, are exactly the flecnodal 
curve, the cuspidal edge and the self-intersection line of the dual surface $S^\vee_\rho$, 
respectively. Since $\F$, $\P$ and $\D$ are parabolas, we state the 

\begin{lemma}\label{parabola_cusp}
The image of the parametrised parabola $t\mapsto (t,ct^2)$, under $\tau_\rho$, 
is the parametrised space curve (lying on $S^\vee$):
$$\a^c_\rho:t\mapsto \left(2(\rho-c)t^3,\ (c-1)t^2,\ 
\left(\frac{c^2}{2}-2c+\frac{3}{2}\rho\right)t^4\right).$$
\end{lemma}

\begin{proof}
This is a direct application of the above Legendre duality map $\tau_\rho$.
\end{proof}

The above parametrisation implies that the curves $\a^c_\rho$ have 
at least $4$-point contact with the $(x,y)$-plane at $t=0$. In order to 
study the behaviour of the curves $\a^c_\rho$ for different values 
of $c$ (for a fixed value of the cr-invariant $\rho$), we will consider 
their projection to the $(x,y)$-plane along the $z$-direction:  
$$\g^c_\rho:t\mapsto \left(2(\rho-c)t^3,\ (c-1)t^2\right).\eqno(6)$$
Clearly, $\g^c_\rho(t)-\a^c_\rho(t)=O(t^4)$. 

\begin{lemma}\label{c=1}
Fix a value of the godron invariant $\rho$. 
The images of all parabolas $y=cx^2$, $c\neq 1$, under the composition of 
$\tau_\rho$ with the projection $(x,y,z)\mapsto (x,y)$, are cusps 
pointing down if $c>1$ and pointing up if $c<1$. These cusps are semi-cubic if 
$c\neq \rho$ and (very) degenerate if $c=\rho$. 

The image of the parabola $y=x^2$ ($c=1$) under the above composition is the 
$x$-axis if $\rho \neq 1$ and it is the origin if $\rho=1$. 
\end{lemma}

\begin{proof}
Lemma~\ref{c=1} and Separating Lemma follow from parametrisation $(6)$.
\end{proof}

\begin{remark*}
It is clear from Lemma~\ref{c=1} that the behaviour of the curve 
$\tau_\rho(\F)$, $\tau_\rho(\P)$ or $\tau_\rho(\D)$ in $S^\vee_\rho$, 
changes drastically when the coefficient $c_{\scriptscriptstyle F}(\rho)$, 
$c_{\scriptscriptstyle P}(\rho)$ or $c_{\scriptscriptstyle D}(\rho)$, 
respectively, passes through the value $1$.
\end{remark*}

\subsection{Proof of Theorem \ref{rho-classification}}
The projection of $S_\rho$ to the $(x,y)$-plane, along the $z$-axis, is a 
local diffeomorphism. So the configuration of the curves $F$, $P$ and $D$ 
with respect to the asymptotic line and the separating $2$-jet at $g$, 
on the surface $S$, is equivalent to the configuration of the parabolas 
$\F$, $\P$ and $\D$ with respect to the parabolas $y=0\cdot x^2=0$ and $y=1\cdot x^2$ 
on the $(x,y)$-plane (see Remark of section~\ref{sec:cr-invariant}). 

\begin{figure}[ht]
\begin{picture}(0,0)
\put(247,55){\bf $\rho$}
\put(164,155){\bf $c$}
\put(121,15){$c_\sD$} 
\put(134,145){$c_\sF$} 
\put(151,0){$c_\sP$} 
\put(93,80){$c_\s=1$}
\end{picture}
\centerline{\psfig{figure=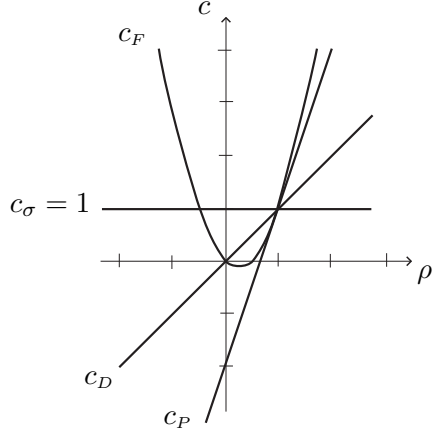,height=5.4cm}}
\caption{\small The coefficients $c_\sF$, $c_\sP$ and $c_\sD$ as functions of the invariant $\rho$.}
\label{coefficients}
\end{figure}

Given a value of $\rho$, this configuration is determined by the order, 
in the real line, of the coefficients of these five parabolas: 
$$c_\sF=\rho(2\rho-1),\ \ c_\sP=(3\rho -2),\ \ c_\sD=\rho,\ \ c_{al}=0,\ \ c_\s=1.$$ 
The graphs of these coefficients, as functions of $\rho$, are depicted in  
Fig.~\ref{coefficients}. 

Using the formulas of the canonical coefficients $c_\sF$, $c_\sP$ and $c_\sD$ 
(or from Fig.~\ref{coefficients}) 
one obtains by straightforward and elementary calculations that : 
$$
\begin{array}{ccc}
\rho\in(1,\infty) &  \iff &  1<c_\sD<c_\sP<c_\sF;\\ 
\rho\in(\frac{2}{3},1) &  \iff &  0<c_\sP<c_\sF<c_\sD<1;\\
\rho\in(\frac{1}{2},\frac{2}{3}) &  \iff &  c_\sP<0<c_\sF<c_\sD<1;\\
\rho\in(0,\frac{1}{2})  & \iff &  c_\sP<c_\sF<0<c_\sD<1;\\
\rho\in(-\frac{1}{2},0) &  \iff &  c_\sP<c_\sD<0<c_\sF<1;\\
\rho\in(-\infty,-\frac{1}{2}) &  \iff &  c_\sP<c_\sD<0<1<c_\sF.
\end{array}
$$
This proves Theorem~\ref{rho-classification}. \fin

\subsection{Proof of Proposition \ref{positive-rho=1}}
Consider the family of surfaces $S_\rho$ given by eq.~$(R)$. 
By eq.~$(P)$, the slope $m$ of the tangent lines of the curve $\P$ 
is given by: $$m(x)=2(3\rho -2)x.$$ 
The slope $p$ of the (double) asymptotic lines along the parabolic curve, 
projected to the $(x,y)$-plane, is given by the equation 
$A^f_p(x,y,p;\frac{\rho}{2})=0$, 
that is, $$p(x)=2x.$$ 

The points of the positive $y$-axis, near the origin, are hyperbolic points 
of the surface $S_\rho$ of eq.~$(R)$. So the hyperbolic domain of $S_\rho$ lies 
locally in the upper side of the parabolic curve. 
Therefore $g$ is a positive (negative) godron if and only if 
the difference of slopes $(p-m)$ is a decreasing (resp. increasing) function of $x$, at $x=0$. 

Consequently, the equation $(p-m)'(0)=-6(\rho -1)$ implies that the godron $g$ 
is positive for $\rho>1$ and negative for $\rho<1$, proving Proposition \ref{positive-rho=1}. \fin

\subsection{Proof of Theorem~\ref{separates} and of Theorem \ref{lr}}
{\bf Preliminary remarks on the asymptotic-double} (see section \ref{index-definition}). \
The asymptotic double $\A$ of the surface $S$ is foliated by the integral 
curves of the asymptotic lifted field of directions. By definition of the lifted 
field, the asymptotic curves of $S$ are the images of these integral curves 
under the natural projection $PT^*S\rightarrow S$ (sending each contact element 
to its point of contact) and, under this projection, the asymptotic double 
$\A$ (of $S$) doubly covers the hyperbolic domain with a fold singularity 
over the parabolic curve. 

Write $\tilde{P}$ for the curve of 
$\A$ which projects over $P$ (that is, the curve formed by the fold points 
in $\A$ of the above projection). 
The surface $\A\setminus \tilde{P}$, has two (not necessarily connected) 
components, noted by $\A_l$ and $\A_r$, separated by $\tilde{P}$. 
The integral curves on the component $\A_l$, 
are projected over the left asymptotic curves and the integral curves on the component 
$\A_r$ are projected over the right ones. 
We call these components the {\em left component} and 
the {\em right component}, respectively, of $\A\setminus \tilde{P}$. 

Now, consider the surface $S$ as the graph of a function 
$f:\R^2 \rightarrow \R$, $z=f(x,y)$, and take the projection 
$\pi: (x,y,z)\rightarrow (x,y)$, along the $z$-axis. The derivative 
of $\pi$ sends the contact elements of $S$ 
onto the contact elements of $\pi(S)\subset \R^2$ and it induces a contactomorphism 
$PT^*S \rightarrow PT^*\R^2$ sending $\A$ to a surface $\tilde\A$ in $PT^*\R^2$, 
which doubly covers (under the natural projection $PT^*\R^2\rightarrow \R^2$) 
the image in $\R^2$ of the hyperbolic domain.  
We still call the surface $\tilde\A\subset PT^*\R^2$ the asymptotic-double of $S$.  
This surface consists of the contact elements of the $(x,y)$-plane 
satisfying the following equation: 
$$f_{xx}dx^2+2f_{xy}dxdy+f_{yy}dy^2=0.\eqno(*)$$
In order to handle the asymptotic double $\tilde\A$, we take an `affine' 
chart of $PT^*\R^2$. The 
space of $1$-jets of the real functions of one real variable $J^1(\R,\R)$ 
(with coordinates $x,y,p$) has a natural contact structure (defined by the 
$1$-form $\a=dy-pdx$) and it parametrises almost all contact elements of 
$\R^2$: The contact element with slope $p_0\neq \infty$ at the point 
$(x_0,y_0)$ of the plane of the variables $(x,y)$ is represented by the 
point $(x_0,y_0,p_0)$ in $J^1(\R,\R)$. The asymptotic-double $\tilde\A$ 
is the surface in $J^1(\R,\R)$ given by the equation 
$$A^f(x,y,p):=f_{xx}+2f_{xy}p+f_{yy}p^2=0,\eqno(7)$$
(obtained from eq.~$(*)$ by taking $p=dy/dx$). 
Moreover, the solutions of the implicit differential equation $(7)$ 
are the images (by $\pi$) of the asymptotic curves of $S$. 
Equation $(7)$ is called the {\em asymptote-equation} of $f$.  

The curve $\tilde{P}$ is the criminant curve (see c.f. \cite{arnoldcs}) 
of the implicit differential equation 
$A^f(x,y,p)=0$ and it is determined by the pair of equations $A^f(x,y,p)=0$ and 
$A^f_p(x,y,p)=0$. 

Below, the images on the plane, under the map $\pi:(x,y,z)\mapsto (x,y)$, 
of the parabolic and flecnodal curves, of the godrons and of the 
hyperbonodes and ellipnodes of $S$, will be called with the same name, 
that is, parabolic curve, etc. One obtains the original objects by applying the 
function $f$ and taking the graph. 
\medskip

\begin{figure}[ht]
\hspace{2cm}{\psfig{figure=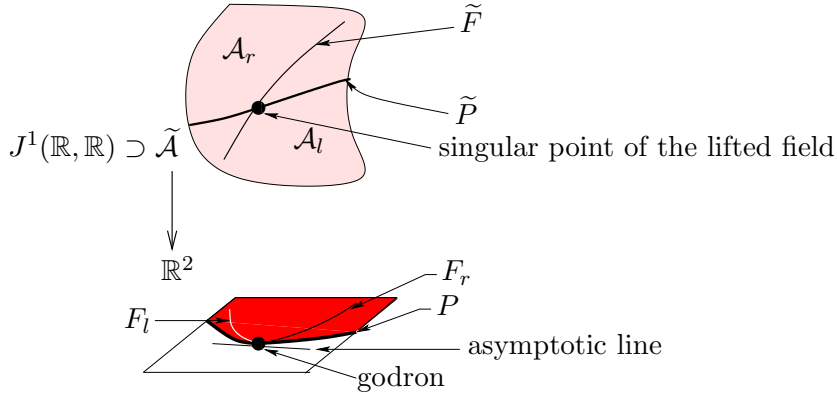,height=5cm}}

\begin{picture}(0,0)
\put(177,145){$\tilde{F}$}
\put(175,110){$\tilde{P}$}
\put(114,101){$\A_l$}
\put(87,135){$\A_r$}
\put(169,98){singular point of the lifted field} 
\put(170,51){$F_r$}
\put(50,33){$F_l$}
\put(168,38){$P$}
\put(180,23){asymptotic line}
\put(138,11){godron} 
\put(6,98){$J^1(\R,\R)\supset\tilde\A$}
\put(64,51){$\R^2$}
\end{picture}

\caption{\small The projection $\pi:\A\rightarrow S$,  
the curves $\tilde{P}$, $\tilde{F}$, $P$ and $F$.}
\label{contour}
\end{figure}

\noindent
{\bf Proof of Theorem~\ref{separates}}.\ ~
Write $\tilde{F}$ for the intersection of $\tilde\A~$ ($A^f(x,y,p)=0$) with the surface given 
by the equation $I^{A^f}(x,y,p)=0$. 
As we mentioned in \S \ref{remarks-computations}, the flecnodal curve in the 
$(x,y)$-plane is the image of the curve $\tilde{F}$ under the 
projection $(x,y,p)\mapsto (x,y)$. The points of (transverse) intersection of the curves 
$\tilde{F}$ and $\tilde{P}$ project over the godrons of $S$. 
So, over a godron the curve $\tilde{P}$ locally separates $\tilde{F}$ 
(see Fig.~\ref{contour}). 

That is, $\tilde{F}$ has one branch on the left component of  $\tilde\A$ and 
other branch on the right component. This implies that 
a godron separates locally the left and right branches of the flecnodal curve, 
proving Theorem~\ref{separates}. \fin
\medskip

\noindent
{\bf Proof of Theorem~\ref{lr}}.\ ~ 
Consider a godron $g$ with cr-invariant $\rho$ of a smooth surface. 
To prove Theorem~\ref{lr} we need to know the values of $\rho$ for which the curves 
$\tilde{P}$ and $\tilde{F}$ are tangent. Of course, such non generic values 
correspond to godrons of non generic surfaces. To found these values we only 
need to know the tangent directions of these curves over $g$. The tangent lines 
of these curves belong to the tangent plane of $\tilde\A$ at $\tilde g$ 
(the point over $g$ in the $(x,y,p)$-space), which is also the contact plane at $\tilde g$. 
So it suffices to take the $4$-jet of $S$ at $g$.

Take the normal form considered above
$$z=\frac{y^2}{2}-x^2y+\rho\frac{x^4}{2}. 
\eqno(R)$$
The coordinates $(x,y,z)$ of this normal form satisfy the conditions considered 
in Theorem~\ref{lr}. 

Since the point $\tilde g$ is the origin, the tangent lines to the curves 
$\tilde{P}$ and $\tilde{F}$ belong to the $(x,p)$-plane. 
The surface $A^f_p(x,y,p)=0$ is the plane given by the equation $p=2x$, which is 
independent of $\rho$. The surface $I^{A^f}(x,y,p)=0$ is the plane given by the 
equation $p=2\rho x$. So the curves $\tilde{P}$ and $\tilde{F}$ are tangent only 
for $\rho=1$ (in this case we have the collapse of two godrons). 

By Proposition~\ref{positive-rho=1}, this implies that the side on which the right 
branch of the flecnodal curve will lie depends only on the index of the godron. 

To see explicitly on which side of the $x$-axis the right branch of the flecnodal 
curve lies for a negative godron, it is enough to look at an example. We will 
take a godron of a cubic surface 
(whose index is $-1$, after Corollary~\ref{cubic}). 

The osculating plane of 
an asymptotic curve at a point of a surface is the tangent plane to 
the surface at that point. 
Using this fact, one defines the ``osculating plane'' of a straight line 
lying in a surface. 

In this way, a segment of a straight line lying in a surface is said to be 
a left (right) curve, if the tangent plane to the surface along that segment 
twists like a left (resp. right) screw.

The $x$-axis is an asymptotic (and flecnodal) curve of the cubic surface 
$z=y^2/2-x^2y$. One verify easily that the positive half axis is a left 
asymptotic curve. This proves Theorem~\ref{lr}. \fin

\subsection{Proof of Theorem \ref{flec-godron_theorem}}

Consider the surface as the graph of a function $z=f(x,y)$. 
In \cite{surf-evolution}, it is proved that the projection of the flecnodal curve to the 
$(x,y)$-plane (the curve $\bar F$ in the above notation) is also the image under the natural projection 
$\pi:\A\subset J^1(\R,\R)\ra\R^2$ ($\pi:(x,y,p)\mapsto (x,y)$) of the 
critical points of the Legendre dual projection 
$\pi^\vee:\A\ni(x,y,p)\mapsto(p,px-y)$ (the {\em folds} of $\pi^\vee$). 
Moreover, it is also proved that, in the hyperbolic domain, the biflecnodes 
correspond to the Whitney pleat singularities of $\pi^\vee:\A\ra(\R^2)^\vee$. 
Both, folds and Whitney pleats, are the 
only stable singularities of a map from a surface to the plane (Whitney).  

It is easy to show that a flec-godron corresponds also to a Whitney pleat of $\pi^\vee$ 
(one shows that, at the point of $\A$ over the flec-godron, the kernel of 
$\pi^\vee_*$ is tangent to the curve of fold points of the map 
$\pi^\vee:\A\ra(\R^2)^\vee$). Thus, after any small generic 
deformation of the surface the flec-godron splits into an ordinary godron and a neighbouring ordinary 
(left or right) biflecnode (since the Whitney pleat singularity is stable). 

The fact that, at the flec-godron moment, both the flecnodal curve and the conodal curve 
have an inflection, follows from the vanishing of the canonical coefficients $c_\sD=\rho$ and 
$c_\sF=\rho(2\rho-1)$, for $\rho=0$. \fin

\subsection{Proof of Theorem \ref{8}}
First, we will prove Theorem~\ref{8} for the case in which the
parabolic curve bounding the hyperbolic disc has only two godrons. 

\begin{lemma}\label{2godrons}
If the parabolic curve bounding a hyperbolic disc $H$ (of a generic smooth 
surface) has exactly two godrons, then the disc $H$ contains an odd number of 
hyperbonodes. 
\end{lemma}

Write $g_1$ and $g_2$ for the godrons lying on $\partial H$. By Proposition~\ref{euler=2},  
both $g_1$ and $g_2$ are positive godrons. 

\begin{claim}\label{claim}
If two vectors $v_1$ and $v_2$ are tangent to $F$ at $g_1$ and $g_2$, respectively, 
and both are pointing from $F_l$ to $F_r$, then $v_1$ and $v_2$ orient the parabolic 
curve $\partial H$ in the same way. 
\end{claim}

\begin{proof}
Since all neighbouring elliptic points of the parabolic curve $\partial H$ belong 
to the same connected component of the elliptic domain, they have the same ``natural'' 
co-orientation (given by the tangent plane). Since both godrons are positive, 
Claim~\ref{claim} follows from Theorem~\ref{lr}. 
\end{proof}

\noindent
{\em Proof of Lemma \ref{2godrons}}.\ ~
Write $f_r$ for the connected component of $F_r$ which starts at $g_1$. 
Since there are only two godrons on $\partial H$, $f_r$ is a segment ending in $g_2$. 
This segment separates $H$ into two parts, which we name $A$ and $B$. 
The connected component of $F_l$ starting in $g_1$, $f_l$, is also a segment ending in $g_2$. 
Claim~\ref{claim} implies 
that if in the neighbourhood of $g_1$ the segment $f_l$ lies in $A$, then, in the neighbourhood 
of $g_2$, it lies in $B$. Thus $f_l$ crosses $f_r$ an odd number of times. 

If $H$ contains other connected components of $F_l$ and $F_r$, then there are (possibly) 
additional hyperbonodes in $H$.
Apart from $f_l$ and $f_r$, the only connected components of $F_l$ and $F_r$ in $H$ are 
closed curves (possibly empty). But the number of intersection points of a closed curve 
of $F_r$ (lying $H$) with $f_l$, or with a closed curve of $F_l$, is even. Thus the number 
of intersection points of $F_l$ with $F_r$ is odd. \fin
\medskip

\noindent
{\bf Proof of Theorem~\ref{8}}.\ ~
To prove the general case of Theorem~\ref{8}, we will consider the closure of the 
hyperbolic disc, the parabolic curve $\partial H$ and the connected components of $F_l$ 
and $F_r$ lying in $H$ as a diagram $\Delta$. We will prove in a purely combinatorial 
manner that the number of intersection points of $F_l$ with $F_r$ is odd. For this, we will 
transform the diagram $\Delta$ using two ``moves'', which are elementary changes 
(of two types) of local diagrams, that preserve the number of intersection points 
of $F_l$ with $F_r$ in the deformed diagram: 

\begin{figure}[ht]
\begin{picture}(0,0)
\put(20,78){(I)} 
\put(19,17){(II)} 
\end{picture}
\centerline{\psfig{figure=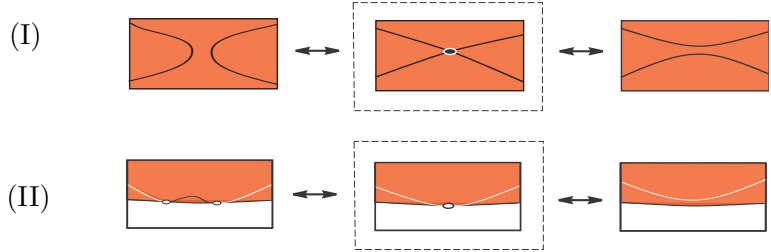,height=3.3cm}}
\caption{\small The two elementary moves of diagrams. 
The moves with opposite choice of colours of the flecnodal curve are also possible.}
\label{moves}
\end{figure}
These moves are depicted in Fig.~\ref{moves}, where an intermediate singular diagram
is marked by a dotted box.

Write $G^+$ and $G^-$ for the number of positive and negative godrons on $\partial H$, 
respectively. Since the asymptotic covering of $\bar{H}$ is a sphere, $G^+-G^-=2$. 

If $G^-=0$, the theorem is proved in Lemma~\ref{2godrons}. So suppose $G^->0$. 

Consider a pair of godrons $g_+$ and $g_-$ of opposite index, 
which are consecutive on $\partial H$. 
Two vectors tangent to $\partial H$ and pointing from $F_l$ to $F_r$, 
one at $g_+$ and the other at $g_-$, provide different orientations of 
$\partial H$ (see Claim~\ref{claim}). 

\begin{figure}[ht]
\begin{picture}(0,0)
\put(59,12){$g_+$}
\put(96,12){$g_-$}
\put(164,13){$g_+$}
\put(196,13){$g_-$} 
\put(267,12){$g_+$}
\put(295,12){$g_-$}
\end{picture}
\centerline{\psfig{figure=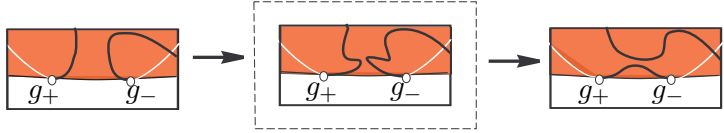,height=1.7cm}}
\caption{\small A deformation of $F_r$ and a move of type I.}
\label{step-1}
\end{figure}

Consider the segment of parabolic curve joining $g_+$ to $g_-$, and which does not 
contain other godrons. 
The local diagram in the tubular neighbourhood of this segment of the parabolic curve 
is depicted in the left side of Fig.~\ref{step-1}. 
\medskip

\noindent
{\bf Step 1}.\ ~
In this tubular neighbourhood we deform the black curves starting in $g_+$ and $g_-$, 
in order to approach one to the other (the central diagram of Fig.~\ref{step-1}). 
Now we apply a move of type I to this diagram in order to obtain a new diagram in which 
the connected component of $F_r$ starting  at $g_+$ will be a segment ending at $g_-$ and 
lying in the tubular neighbourhood of the considered segment of the parabolic curve. 

\noindent
(It must be clear that, in Fig.~\ref{moves} and in Fig.~\ref{step-1}, we are not claiming 
that the surface is being deformed in such a way that this evolution of 
the flecnodal curve happens. We are deforming the diagram, not the surface. 
However, to avoid new notations and symbols, we have kept the names: godron, 
parabolic curve, etc; and the notations:  $g_+$, $g_-$, $F_l$, $F_r$, $\partial H$, etc.)
\medskip

\noindent
{\bf Step 2}.\ ~Applying a move of type II to the local diagram obtained in Step~1, 
one obtains a new diagram without the pair of godrons $g_+$ and $g_-$. 
\medskip

Applying $G^-$ times the above process, one obtains a final diagram having only two positive godrons.  
Theorem~\ref{8} is proved applying Lemma~\ref{2godrons} to this final diagram (note that also the 
proof of Lemma~\ref{2godrons} depends only on the combinatorial properties of the 
initial diagram). \fin

\subsection{Proof of Theorem \ref{4-fronts}}

To prove Theorem~\ref{4-fronts} we will use the fact that 
{\em the tangential map of $S$ sends the 
flecnodal curve of $S$ onto the flecnodal curve of $S^\vee$}. 

The dual of a front $\tilde{S}$ in general position at a swallowtail point $s$ is a godron 
of a (locally) smooth surface. So Theorem~\ref{tangency} and the Separating 
Lemma imply that the flecnodal curve of $\tilde{S}$ has a cusp at $s$ having the same 
tangent line that the cuspidal edge of $\tilde{S}$. Now, by Theorem~\ref{separates}, 
the swallowtail point 
separates the flecnodal curve into its left and right branches. 

The configuration formed by the flecnodal curve, the cuspidal edge and the self-intersection 
line of $\tilde{S}$ at the swallowtail point $s$, is determined by the configuration 
formed by the curves $F$, $P$, $D$ and the separating $2$-jet on the (locally smooth) dual surface 
$\tilde{S}^\vee$, at its godron $g=s^\vee$. 

Theorem~\ref{rho-classification} says that, at a godron $g$, there are six possible generic 
configurations of the curves $F$, $P$ and $D$, with respect to the separating $2$-jet 
and to the asymptotic line at $g$. 
Since the asymptotic line is not considered
in the concerned configurations, we can eliminate the number $0$ (corresponding to the 
asymptotic line) from the six inequalities of the proof of 
Theorem~\ref{rho-classification}. One obtains four distinct inequalities, 
corresponding to four open intervals for the values of $\rho$ : 

$$
\begin{array}{lcc}
\rho\in(1,\infty) &  \iff &  1<c_\sD<c_\sP<c_\sF;\\ 
\rho\in(0,1)  & \iff &  c_\sP<c_\sF<c_\sD<1;\\
\rho\in(-\frac{1}{2},0) &  \iff &  c_\sP<c_\sD<c_\sF<1;\\
\rho\in(-\infty,-\frac{1}{2}) &  \iff &  c_\sP<c_\sD<1<c_\sF.
\end{array}
$$
Using the Separating Lemma and the configurations 
of Theorem~\ref{rho-classification} (not considering the asymptotic line) one obtains 
that these four configurations correspond to the four configurations (of 
Theorem~\ref{4-fronts}) for the flecnodal curve, the cuspidal edge and the 
self-intersection line in the neighbourhood of a swallowtail 
point of a front in general position. \fin

\subsection{Proof of Theorems~\ref{7-fronts}, \ref{4/3} and 
\ref{rho-classification_contour}}\label{proofs_4/3}

The following general lemma is the key to prove of Theorems~\ref{7-fronts}, \ref{4/3} 
and  \ref{6-fronts}. 

Let $q$ be a point of a front $S$ of $\RP^3$, in general position. 
Write $q^\vee$ for the point of the dual surface $S^\vee\subset(\RP^3)^\vee$, 
corresponding to $q$ (that is, $q^\vee$ is the tangent plane to $S$ at $q$), 
and write $\Pi$ for the plane tangent to $S^\vee$ at $q^\vee$. 

\begin{qlemma}
The image of the (total) $q$-contour of $S$ under the tangential map 
$S\ra S^\vee\subset (\RP^3)^\vee$ is the tangent section $\Pi\cap S^\vee$ 
of the dual surface $S^\vee$ at $q^\vee$. 
\end{qlemma}

\begin{proof}
One needs to prove that any plane tangent to $S$ and passing through $q$ is a 
point of $\RP^3$ belonging to $\Pi\cap S^\vee$, and that every point of $\Pi\cap S^\vee$ 
is a plane tangent to $S$ passing through $q$. 

Since the point $q\in S\subset \RP^3$ is precisely the plane 
$\Pi$ of $(\RP^3)^\vee$ tangent to $S^\vee$ at $q^\vee$, the points 
of $\Pi$ are the planes of $\RP^3$ passing through $q$. Consequently, the 
points of $(\RP^3)^\vee$ belonging to $\Pi\cap S^\vee$ 
(forming the tangent section of $S^\vee$ at $q^\vee$) 
are the planes of $\RP^3$ tangent to $S$ (since they belong to $S^\vee$) and 
passing through $q$ (since they belong to $\Pi$). 
\end{proof}

The $q$-Contour Lemma and the arguments given in the proof of Theorem~\ref{4-fronts} 
imply that Theorems~\ref{7-fronts} and \ref{4/3} together 
are equivalent to Theorem~\ref{rho-classification_contour}. So we only need to 
prove Theorem~\ref{rho-classification_contour}. 
\medskip

\noindent
{\bf Proof of Theorem~\ref{rho-classification_contour}}. \ 
Consider the surface $S_\rho$ given by eq.~$(R)$. We will use the parametrisation 
$(R^\vee)$ of the dual surface $S_\rho^\vee$ given in \S\ref{proof_separating}. 

So, in order to find the first terms of the power series expansion of the $g$-contour 
of $S$,  we need to find the zeros of the equation $y^2/2-2x^2y+3\rho x^4=0$. Completing squares, 
one easily factorises the left hand side of this equation: 
$$\frac{1}{2}\left(y-(2+\sqrt{4-3\rho})x^2\right)\left(y-(2-\sqrt{4-3\rho})x^2\right)=0.$$

This implies that the projection of the $g$-contour of $S$ to the $(x,y)$-plane 
(in the $z$-direction) is given by two tangent curves: 
$$C_-: \ y=c_\sC^-x^2+\ldots \ \ \ \mbox{ and }\ \ \ \ C_+: \ y=c_\sC^+x^2+\ldots,$$ 
where $~c_\sC^-=2-\sqrt{4-3\rho}~$ and $~c_\sC^+=2+\sqrt{4-3\rho}$.

\begin{figure}[ht]
\begin{picture}(0,0)
\put(167,55){\bf $\rho$}
\put(42,167){$c_\sC^+$}
\put(23,142){$\scriptstyle{\rho=-\frac{\sqrt{7}}{2}}$}
\put(80,168){\bf $c$}
\put(42,40){$c_\sC^-$}
\put(37,18){$c_\sD$} 
\put(57,120){$c_\sF$} 
\put(150,137){$\scriptstyle{\rho=\frac{\sqrt{7}}{2}}$} 
\put(149,105){$\scriptstyle{\rho=\frac{4}{3}}$} 
\put(69,0){$c_\sP$} 
\put(7,84){$c_\s=1$}
\put(99,-4){$\bf C$}

\put(330,55){\bf $\rho$}
\put(246,168){\bf $c$}
\put(201,120){$c_\sT^+$}
\put(197,95){$\scriptstyle{\rho=-\frac{\sqrt{3}}{2}}$}
\put(201,44){$c_\sT^-$}
\put(201,18){$c_\sD$} 
\put(215,150){$c_\sF$} 
\put(233,0){$c_\sP$} 
\put(320,84){$c_\s=1$}
\put(303,74){$\scriptstyle{\rho=\frac{8}{9}}$} 
\put(295,50){$\scriptstyle{\rho=\frac{\sqrt{3}}{2}}$}
\put(266,-4){$\bf T$}
\end{picture}
\centerline{\psfig{figure=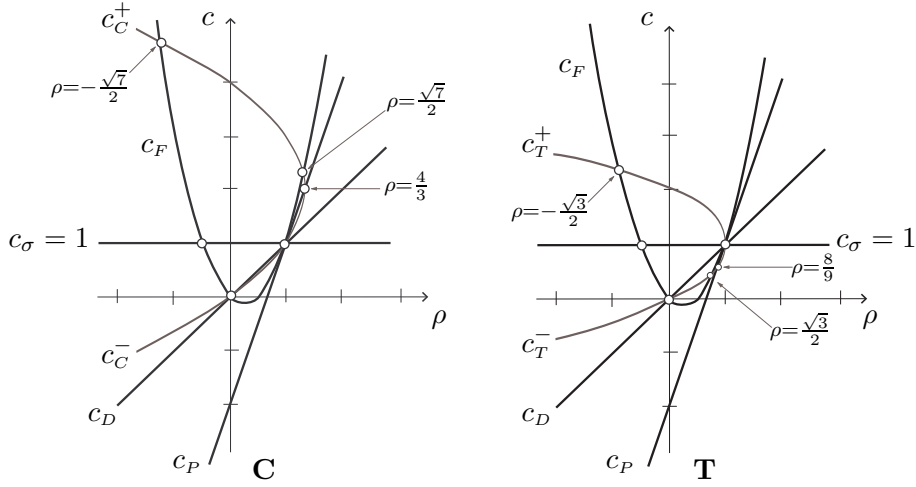,height=6cm}}
\caption{\small The coefficients $c_\sF$, $c_\sP$ and $c_\sD$ as functions of the invariant $\rho$, 
together with the coefficients $c_\sC^-$, $c_\sC^+$ (part C) and $c_\sT^-$, $c_\sT^+$ (part T).}
\label{coefficients_C_T}
\end{figure}

Of course these two coefficients are defined as functions of $\rho$ only 
for $\rho<\frac{4}{3}$. The graphs of these coefficients, as functions of $\rho$, 
are depicted in Fig.~\ref{coefficients_C_T}.C (together with the graphs of $c_\sF$, 
$c_\sP$, $c_\sD$ and $c_\s=1$). Indeed, the union of these two graphs 
forms the parabola given by the equation:
$$\rho=-\frac{1}{3}(c-2)^2+\frac{4}{3}.$$ 

Using the above formulas of $c_\sC^-$ and $c_\sC^+$, and those of the canonical coefficients 
$c_\sF$, $c_\sP$, $c_\sD$ and $c_\s$: 
$$c_\sF=\rho(2\rho-1),\ \ c_\sP=(3\rho -2),\ \ c_\sD=\rho,\ \ c_\s=1,$$ 
(or from Fig.~\ref{coefficients_C_T}.C) 
one obtains the inequalities of Theorem~\ref{rho-classification_contour} (and the corresponding 
exceptional values of $\rho$: $\frac{4}{3}$, $\pm\frac{\sqrt{7}}{2}$, $-\frac{1}{2}$, $0$, $1$) 
by straightforward and elementary calculations. To make the calculations simpler, one can use 
the inequalities already obtained in Theorem~\ref{rho-classification} (eliminating the number $0$). \fin
\medskip

\begin{remark*}
One can also prove Theorem~\ref{rho-classification_contour} by using the parametrisations 
(obtained from Lemma~\ref{parabola_cusp} of \S\ref{proof_separating}) of the 
self-intersection line $D^\vee$, the flecnodal curve $F^\vee$ and the cuspidal edge 
$\Sigma=P^\vee$ of the dual swallowtail $S^\vee$:  
\smallskip

$D^\vee=(0,(\rho-1)t^2, \frac{1}{2}\rho(\rho-1)t^4),$ 

$F^\vee=(4\rho(1-\rho)t^3, (\rho(2\rho-1)-1)t^2, \frac{1}{2}\rho(\rho-1)(4\rho^2-7)t^4),$ 

$P^\vee=(4(1-\rho)t^3, 3(\rho-1)t^2, 3(\rho-1)(3\rho-4)t^4).$
\smallskip

For one needs to find the values of $\rho$ for which a component of one of these curves 
change its sign, and then to use the geometric meaning of that sign change. 

For instance, 
the flecnodal curve of the swallowtail $S^\vee$ has higher order of contact with the 
tangent plane (passing locally from one side to the other of it) when the third 
component of $F^\vee$, $~\frac{1}{2}\rho(\rho-1)(4\rho^2-7)t^4,~$ equals zero. 
In this case one finds the degenerate case $\rho=1$, the flec-godron $\rho=0$ 
and the exceptional values $\pm\frac{\sqrt{7}}{2}$. 
\end{remark*}

\subsection{Proof of Theorems~\ref{6-fronts} and \ref{rho-classification_tangent}}

As in \S\ref{proofs_4/3}, the $q$-Contour Lemma and the arguments given 
in the proof of Theorem~\ref{4-fronts} imply that Lemma~\ref{tangent_section}
and Theorem~\ref{rho-classification_tangent} together are equivalent to Theorem~\ref{6-fronts}. 
So, we will prove only Lemma~\ref{tangent_section} and Theorem~\ref{rho-classification_tangent}. 
\medskip

\noindent
{\bf Proof of Lemma~\ref{tangent_section} and Theorem~\ref{rho-classification_tangent}}. \ 
Consider the surface $S_\rho$ given by the equation $z=y^2/2-x^2y+\rho x^4/2+\f(x,y)$, 
where $\f(x,y)$ is the sum of monomials in $x$ and $y$ of degree greater than $4$.  
The first terms of the power 
series expansion of the branches of the tangent section at $g$ are thus
given by the zeros of the equation $y^2/2-x^2y+\rho x^4/2=0$. Completing squares, 
we factorise the left hand side of this equation:
$$\frac{1}{2}\left(y-(1+\sqrt{1-\rho})x^2\right)\left(y-(1-\sqrt{1-\rho})x^2\right)=0.$$
Hence, the tangent section of $S$ at $g$ is given by two (simply) tangent curves: 
$$C_-: \ y=c_\sT^-x^2+\ldots \ \ \ \mbox{ and }\ \ \ \ C_+: \ y=c_\sT^+x^2+\ldots,$$ 
where $~c_\sT^-=1-\sqrt{1-\rho}~$ and $~c_\sT^+=1+\sqrt{1-\rho}~$ (proving Lemma~\ref{tangent_section}).

The coefficients $c_\sT^-$ and $c_\sT^+$ are defined as functions of $\rho$ only 
for $\rho<1$, and their graphs, as functions of $\rho$, 
are depicted in Fig.~\ref{coefficients_C_T}.T (together with the graphs of $c_\sF$, 
$c_\sP$, $c_\sD$ and $c_\s=1$). The union of these two graphs 
forms the parabola given by the equation:
$$\rho=-(c-1)^2+1.$$ 

Using the above formulas of $c_\sT^-$ and $c_\sT^+$, and those of the coefficients 
$c_\sF$, $c_\sP$, $c_\sD$ and $c_\s$ (or from Fig.~\ref{coefficients_C_T}.T) 
one obtains the inequalities of Theorem~\ref{rho-classification_tangent} (and the corresponding 
exceptional values of $\rho$: $\frac{8}{9}$, $\pm\frac{\sqrt{3}}{2}$, $-\frac{1}{2}$, $0$, $1$) 
by straightforward and elementary calculations. \fin
\medskip

\begin{remark*}
One can also prove Theorem~\ref{rho-classification_tangent} by using the parametrisations 
of the conodal curve $D$, the flecnodal curve $F$ and the parabolic curve
$P$ of $S$:  
\smallskip

$D=(x,\rho x^2+\ldots, \frac{1}{2}\rho(\rho-1)x^4+\ldots),$ 

$F=(x ,\rho(2\rho-1)x^2+\ldots,  \frac{1}{2}\rho(\rho-1)(4\rho^2-3)x^4+\ldots),$ 

$P=(x, (3\rho-2)x^2+\ldots, \frac{1}{2}\rho(\rho-1)(9\rho-8)x^4+\ldots)$. 
\smallskip

For one needs to find the values of $\rho$ for which the first term of 
a component of one of these curves change its sign, and then to use the geometric 
meaning of that sign change. 

For example, 
the flecnodal curve of $S$ has higher order of contact with the 
tangent plane (passing locally from one side to the other of it) when the 
term $~\frac{1}{2}\rho(\rho-1)(4\rho^2-3)t^4,~$ of the third 
component of $F$, equals zero, providing the degenerate case $\rho=1$, 
the flec-godron $\rho=0$  and the exceptional values 
$\pm\frac{\sqrt{3}}{2}$. 
\end{remark*}

\begin{remark*}
When this paper was almost finished, I visited l'\'Ecole Normale Sup\'erieure de Lyon 
to give a talk about the results of \cite{surf-evolution} and of this paper. Few days before 
my talk, E.~Ghys and D.~Serre have found the book \cite{Levelt} on the history of 
thermodynamics in Netherlands. It describes a part of Korteweg's work
(\cite{Kortewegpp,Korteweggtp}) about the godrons (called plaits in \cite{Levelt}), 
the parabolic curve and the conodal curve. According to \cite{Levelt}, Korteweg had also 
described the bifurcations of the parabolic and conodal curves when two godrons 
are born or disappear, for an evolving 
surface. 
Korteweg's work on the theory of surfaces was motivated by 
thermodynamical problems. 
\end{remark*}


\end{document}